\documentclass[11pt]{amsart}
\usepackage{amscd,amsmath,amssymb,amsfonts,verbatim}
\usepackage[cmtip, all]{xy}
\usepackage{comment}
\usepackage{tikz-cd}
\usepackage{hyperref} 
\usepackage{float}  
\usepackage{dynkin-diagrams}
\usepackage{geometry}
\renewcommand\qedsymbol{$\blacksquare$}

\newtheorem{thm}{Theorem}[section]
\newtheorem{prop}[thm]{Proposition}
\newtheorem{lem}[thm]{Lemma}
\newtheorem{cor}[thm]{Corollary}

\theoremstyle{definition}

\newtheorem{defn}[thm]{Definition}
\newtheorem{remks}[thm]{Remarks}
\newtheorem{remk}[thm]{Remark}

\newtheorem{exm}[thm]{Example}
\newtheorem{exms}[thm]{Examples}
\newtheorem{notat}[thm]{Notation}
\numberwithin{equation}{section}
\numberwithin{equation}{section}
\newtheorem{ques}[thm]{Question}

{\hfill$\square$\end{defn}}
\newenvironment{rem}{\begin{remk}}%
{\hfill$\square$\end{remk}}
{\hfill$\square$\end{remks}}
{\hfill$\square$\end{exm}}
{\hfill$\square$\end{exms}}
{\hfill$\square$\end{notat}}

\newcommand{\Z}{{\mathbb Z}}

\newcommand{\rank}{{\rm rank}}

\newcommand{\Spec}{{\rm Spec}}

\newcommand{\Ker}{{\rm Ker}}
\newcommand{\Img}{{\rm Im}}

\newcommand{\hh}{{\rm ht}}

\newcommand{\Um}{{\rm Um}}

\newcommand{\Aut}{{\rm Aut}}

\newcommand{\SR}{{\rm SR}}
\newcommand{\sr}{{\rm s.r.}}

\renewcommand{\>}{\rangle}

\newcommand{\SL}{{\operatorname{\rm SL}}}

\newcommand{\Spin}{{\operatorname{\rm Spin}}}
\newcommand{\ESpin}{{\operatorname{\rm ESpin}}}
\newcommand{\St}{\mathrm{St}}

\newcommand{\KSpin}{{\operatorname{\rm KSpin}}}
\newcommand{\KO}{{\operatorname{\rm KO}}}
\newcommand{\K}{{\operatorname{\rm K}}}
\newcommand{\SK}{{\operatorname{\rm SK}}}

\newcommand{\A}{{\operatorname{\rm A}}}
\newcommand{\B}{{\operatorname{\rm B}}}
\newcommand{\C}{{\operatorname{\rm C}}}
\newcommand{\D}{{\operatorname{\rm D}}}
\newcommand{\E}{{\operatorname{\rm E}}}
\newcommand{\F}{{\operatorname{\rm F}}}
\newcommand{\G}{{\operatorname{\rm G}}}
\newcommand{\Or}{{\operatorname{\rm O}}}

\newcommand{\ds}{{/\kern-3pt/}}

\renewcommand{\dim}{\text{\rm dim}}

\newcommand{\tuborg}{\left\{\begin{array}{ll}}
\newcommand{\sluttuborg}{\end{array}\right.}

\newcounter{elno}

\newcounter{elno-abc}   

\newcounter{elno-abc-prime}

\begin{document}

\author{Chandan Bhaumik}
%\title[Stability theorems for simply-laced $\K_{2}$ functor]{Stability theorems for simply-laced $\K_{2}$ functor over Laurent polynomial rings}
\title[Stability for $\K_{2}$ functor of simply-laced Chevalley groups]{Stability for $\K_{2}$ functor of simply-laced Chevalley groups over Laurent polynomial rings}
%\author{Chandan Bhaumik}
%\author{Chandan Bhaumik\textsuperscript{*} and Manoj K. Keshari}

\address{Department of Mathematics, Indian Institute of Technology Bombay, Powai, Mumbai 400076, India}
\email{cbhaumik11math@gmail.com, cbhaumik@math.iitb.ac.in}

% \address{(H.P. Sarwar) Department  of Mathematics, Indian  Institute of Technology Kharagpur, Kharagpur 721302, West Bengal, India}
% \email{mathparvez@gmail.com, parvez@maths.iitkgp.ac.in}

%\thanks{\textsuperscript{*}Corresponding author. Email: cbhaumik11math@gmail.com (Chandan Bhaumik)}

\dedicatory{}



\maketitle
% \vspace{-.7cm}
% \begin{center}
%    \small Department of Mathematics, IIT Bombay, Powai, Mumbai-400076, India\\ cbhaumik11math@gmail.com
% \end{center}
% \vspace{.5cm}
%\noindent\textit{\textsuperscript{*}Corresponding author. Email: cbhaumik@math.iitb.ac.in}

\begin{quote}\emph{Abstract.} 
This article is devoted to establishing stability theorems for the $\K_{2}$-functor of Chevalley groups $\G(\Phi,-)$ of simply-laced type root systems $\Phi$ over Laurent polynomial extensions of a Noetherian domain. As an application, we obtain a $\K_2$-homotopy invariance property for Laurent polynomial rings.
%This paper is devoted to the study of the stability theorems for the $\K_{2}$-functor of Chevalley groups $\G(\Phi,-)$ of simply-laced type root systems $\Phi$ over Laurent polynomial extensions of a Noetherian domain.
\end{quote}
%\end{abstract}
\setcounter{tocdepth}{1} % Ensures subsections appear in the table of contents

%\tableofcontents

\section{Introduction}\label{sec:Intro}
The stabilization problem has been of great interest throughout the development of algebraic $\K$-theory, initiated by Bass \cite{Ba68}. This problem encompasses two aspects, namely the surjectivity and injectivity of the stabilization maps. Bass proposed that the stability conditions should be formulated in terms of the relationship between the dimension of the maximal spectrum of a ring and the rank of a maximal split torus of the underlying algebraic group. This criterion is optimal, as evidenced by counterexamples to non-stability for the $\K_{2}$-functor arising from the homotopy theory of real orthogonal groups (see \cite[Section 7]{Ka76}). A primary objective of the present article is to investigate the stabilization problem for the $\K_{2}$-functor of Chevalley groups associated with simply-laced type root systems over Laurent polynomial extensions of a Noetherian domain under the same stability condition. 
\par
First, we recall the notion of $\K_{2}$-stability for Chevalley groups. Let $\Phi$ be a reduced irreducible root system, and let $\G(\Phi,-)$ be a simply connected Chevalley Demazure group scheme over $\Z$ of type $\Phi$ in the sense of \cite{DeGr70}. If $R$ is a commutative ring with $1$, the value of the functor $\G(\Phi,-)$ on $R$ is called the simply connected Chevalley group of type $\Phi$ over $R$, and is denoted by $\G(\Phi,R)$. The subgroup of $\G(\Phi,R)$ generated by the root unipotents $e_{\alpha}(r)$ for all $\alpha\in\Phi$ and all $r\in R$ is called the elementary subgroup of $\G(\Phi,R)$, and is denoted by $\E(\Phi,R)$. Let $\St(\Phi,R)$ be a Steinberg group presented by generators $x_{\alpha}(r)$ (see Definition $\ref{defn_stein_grp}$); there is a natural surjective homomorphism $\St(\Phi, R) \twoheadrightarrow \E(\Phi, R)$ defined by $x_{\alpha}(r)\mapsto e_{\alpha}(r)$. Its kernel is defined to be the $\K$-theoretic group $\K_2(\Phi, R)$, which is central in $\St(\Phi, R)$ when $\rank(\Phi)\geq 3$ \cite{LaSiVo24}.\par
 Any (closed) embedding of root systems $\Phi'\subseteq \Phi$ induces a homomorphism of Steinberg groups $\St(\Phi',R)\to\St(\Phi,R)$, and hence a homomorphism of the corresponding $\K_{2}$-functors $\K_{2}(\Phi',R)\to\K_{2}(\Phi,R)$. The surjective (resp. injective) stability problem is to determine conditions on the ring $R$ depending on the embedding $\Phi'\subseteq \Phi$, under which the stabilization map $\K_{2}(\Phi',R)\to\K_{2}(\Phi,R)$ is surjective (resp. injective). The $\K_2$-stability provides criteria for determining when the unstable group $\K_2(\Phi,R)$ has stabilized as the rank of the root system increases. For a classical root system $\Phi$, such stability results connect the unstable group $\K_2(\Phi,R)$ associated with a root system of finite rank to the stable group $\K_2(\Phi_{\infty},R)$. \par
 Numerous $\K_{2}$-stability results are available in the literature. For Milnor's $\K_{2}$-groups, surjective stability was first established by Dennis \cite{Da73} and Vaserstein \cite{Va75}, while injective stability was proved by van der Kallen \cite{Ka76} and Suslin--Tulenbaev \cite{SuTu76}.  In the broader context of Chevalley groups, this problem was studied by Stein \cite{St78}, who provided a versatile method that handles both classical and exceptional types uniformly. Tulenbaev \cite{Tu83} improved the stability result for Milnor $\K_{2}$-groups over polynomial rings in several variables.  \par
 The first main contribution of this paper is an extension of Tulenbaev's stability result for Milnor $\K_{2}$-groups to Chevalley $\K_{2}$-groups associated with simply-laced type root systems. This is captured in the following theorem.
% which generalizes \cite[Theorem 5.3]{Tu83} and \cite[Theorem 5.9]{Si25}. 
\begin{thm}[Theorems \ref{thm_surj_stab_poly}, \ref{thm_inj_stab_poly}]\label{intro_thm_stab_poly}
    Let $A=R[X_{1},\ldots,X_{n}]$ be a polynomial extension of a commutative Noetherian domain $R$. Suppose $\Phi'\subset \Phi$ is one of the following embeddings of reduced irreducible simply-laced root systems:
    \begin{enumerate}
        \item[$(i)$] $\A_{r}\subset \A_{r+1}$, with $r\geq\max\{3,\dim(R)+1\}$;
        \item[$(ii)$] $\D_{r}\subset \D_{r+1}$, with $r\geq\max\{4,\dim(R)+2\}$;
        \item[$(iii)$] $\E_{r}\subset \E_{r+1}~(r=5,6,7)$, with $\dim(R)\leq r-3$;
       \item[$(iv)$] $\D_{r}\subset \E_{r+1}~(r=6,7)$ and $\D_{8}\subset \E_{8}$, with $\dim(R)\leq 2$;
        \item[$(v)$] $\A_{r}\subset \E_{r+1}~(r=5,6,7)$ and $\A_{r}\subset \D_{r+1}~(r\geq 5)$, with $\dim(R)\leq 1$;
         \item[$(vi)$] $\A_{r}\subset \E_{r}~(r=7,8)$, with $\dim(R)\leq 1$;
          \item[$(vii)$] $\A_{4}\subset\A_{\geq 4},~\A_{4}\subset\D_{\geq 5}$ and $\A_{4}\subset\E_{6,7,8}$, with $\dim(R)\leq 1$;
          \item[$(viii)$] $\A_{r}\subset \A_{r+1}$, with $r\geq\max\{4,\dim(R)+2\}$;
        \item[$(ix)$] $\D_{r}\subset \D_{r+1}$, with $r\geq 5$.
    \end{enumerate}
    In $(ix)$, assume further that $R$ is either a principal ideal domain or a Dedekind domain containing a field. Then, under the given conditions, the canonical homomorphism $\K_{2}(\Phi',A)\to\K_{2}(\Phi,A)$ is surjective in cases $(i)-(vii)$ and injective in cases $(viii)-(ix)$. 
\end{thm}
%In the case $n=0$, parts (i)--(iii) of the above theorem were obtained by Stein \cite{St78}.
For arbitrary $n$, the injective and surjective stability results (i) and (vii) are due to Tulenbaev \cite[Theorem 5.3]{Tu83}, and the surjective stability result (vii) is due to Sinchuk \cite[Theorem 5.9]{Si25}. However, the remaining cases covered by Theorem \ref{intro_thm_stab_poly} were not previously known and are therefore genuinely new.\par 
It is natural to ask whether one can extend the above theorem to the setting of Laurent polynomial extensions. It is worth noting that, following Tulenbaev's contributions, no further progress on stability questions over Laurent polynomial rings has been reported. In particular, neither injective nor surjective stability has been established in this context for the $\K_{2}$-functor of Chevalley groups. Motivated by this gap, our second main contribution is the following result.
\begin{thm}[Theorems \ref{thm_surj_stab_Laur_poly}, \ref{thm_inj_stab_Laur_poly}]\label{intro_thm_stab_Laur_poly}
    Let $A=R[X_1,\ldots,X_n,Y_{1}^{\pm1},\ldots,Y_{m}^{\pm1}]$ be a Laurent polynomial extension of a commutative Noetherian domain $R$, and let $\Phi'\subset \Phi$ be embeddings of root systems as in Theorem $\ref{intro_thm_stab_poly}$. Consider the canonical homomorphism $\Gamma_A:\K_{2}(\Phi',A)\to\K_{2}(\Phi,A)$. Then the following holds: \begin{enumerate}
        \item[(1)] If $n=0,m=1$ and assume that $\K_{1}(\Phi,R[Y_{1}])$ is trivial, then $\Gamma_{A}$ is surjective in cases $(i)-(vii)$.
        \item[(2)] $\Gamma_{A}$ is injective in case $(viii)$.
    \end{enumerate} 
\end{thm}
As an application of the stability theorems established above, we obtain the following unstable $\K_2$-homotopy invariance result for Laurent polynomial rings.
\begin{thm}[Theorem \ref{thm_hom_inv_Laur_poly}]\label{intro_thm_hom_inv_Laur_poly}     Let $R$ be a Dedekind domain or field, having no non-trivial units. Let $\Psi$ be a root system of type $\A_{\geq 4},~\D_{\geq 5}$ or $\E_{6,7,8}$. Assume that $\K_{1}(\Psi,R)$ is trivial. Then $\K_{2}(\Psi,R[X,X^{-1}])\cong \K_{2}(\Psi,R)$. Moreover, if $\K_{2}(\Psi,R)$ is generated by Steinberg symbols (resp. Dennis-Stein symbols), then so is $\K_{2}(\Psi,R[X,X^{-1}])$.
\end{thm}
When $\Psi$ is an irreducible simply-laced root system containing $\A_4$ but not of type $\E_8$, the unstable $\K_2$-homotopy invariance result for polynomial extensions of a Dedekind domain is due to Sinchuk \cite[Theorem 1.3]{Si25}. It is noted that if $R$ is a field and $\Phi$ is a reduced irreducible root system, then $\K_{2}(\Phi,R[X,X^{-1}])$ is generated by Steinberg symbols (and hence Dennis-Stein symbols) due to Abe-Morita \cite[Corollary 6]{AbMo88}. 
%For Milnor's $\K_{2}$-groups, the stable $\K_2$-homotopy invariance result for $c$-divisible monoids without nontrivial units over a regular ring is due to Mushkudiani \cite{Mu95}. For $\Phi$-simplicial $c$-divisible monoids, the corresponding result is due to Gubeladze \cite{Gu90}. \par
%The remainder of this article is organized as follows. In Section $\ref{sec_preli}$, we collect the basic definitions and preliminary results that will be used throughout the article. In Section $\ref{sec_mip}$, we establish the monic inversion principles that play a crucial role in the proofs of our main theorems. Using these tools, Section $\ref{sec_pr}$ addresses the stabilization problem and proves Theorem $\ref{intro_thm_stab_poly}$. Finally, in Section $\ref{sec_lpr}$, we extend our analysis to Laurent polynomial extensions and prove Theorems $\ref{intro_thm_stab_Laur_poly}$ and $\ref{intro_thm_hom_inv_Laur_poly}{}$. 
\section{Steinberg Groups and $\K_{2}$ Functor}\label{sec_preli}
{\bf Throughout the article, all rings are assumed to be commutative with unity, and root systems are taken to be reduced and irreducible.}
\begin{defn}
\label{defn_stein_grp}
Let R be a ring and $\Phi$ be a root system. The {\it Steinberg group} $\St(\Phi, R)$ is the group generated by symbols $x_\alpha(u)$, where $\alpha \in \Phi$ and $u \in R$, subject to the following relations:
\begin{align}
    &    x_\alpha(u) \, x_\alpha(v) = x_\alpha(u+v), \qquad u, v \in R,~ \alpha\in\Phi;\label{St1}\\[4pt]
    &    [x_\alpha(u), x_\beta(v)] = \prod_{\substack{i,j > 0 \\ i\alpha + j\beta \in \Phi}} x_{i\alpha + j\beta}(C_{ij\alpha\beta}\, u^i v^j), \qquad u, v \in R,~\alpha, \beta \in \Phi, ~\alpha\neq -\beta,\label{St2}
\end{align}
where $[a,b] = aba^{-1}b^{-1}$, the product is taken over all positive integers $i,j$ such that $i\alpha + j\beta \in \Phi$, and $C_{ij\alpha\beta}\in\Z$. The constants $C_{ij\alpha\beta}$ (depending only on $\Phi$) are called the Chevalley structure constants. 
\end{defn} 
For a ring $R$ and a root system $\Phi$, let $\G(\Phi,R)$ denote the simply connected {\it Chevalley group} of type $\Phi$, and $\E(\Phi,R)$ denote its elementary subgroup generated by elementary root unipotents $e_{\alpha}(r)$, $\alpha\in\Phi,~r\in R$. Generators $e_{\alpha}(r)\in \E(\Phi,R)$ satisfy the relations $(\ref{St1})-(\ref{St2})$, whereas $x_{\alpha}(r)$ are regarded as formal symbols. There is a well-defined map $\eta^{\Phi}:\St(\Phi,R)\to\G(\Phi,R)$ sending $x_{\alpha}(r)\mapsto e_{\alpha}(r)$, and the image of $\eta^{\Phi}$ is $\E(\Phi,R)$, which is a normal subgroup of $\G(\Phi,R)$ when $\rank(\Phi)\geq 2$. The groups $\K_{1}(\Phi,R)$ and $\K_{2}(\Phi,R)$ are defined as the cokernel and kernel of $\eta^{\Phi}$, respectively (see \cite{St78})
\begin{center}
        \begin{tikzcd}
1 \arrow[r] & \K_{2}(\Phi,R) \arrow[r]           & \St(\Phi, R) \arrow[r,"\eta^{\Phi}"]                 & \G(\Phi, R) \arrow[r]  & \K_{1}(\Phi,R) \arrow[r]               & 1. 
\end{tikzcd}
\end{center}
When $\Phi=\A_r$, the Chevalley group $\G(\Phi,R)$ is the special linear group $\SL_{r+1}(R)$, and its elementary subgroup $\E(\A_r,R)$ coincides with the linear elementary group $\E_{r+1}(R)$. While $\G(\D_r,R)=\Spin_{2r}(R)$ and its elementary subgroup $\E(\D_r,R)=\ESpin_{2r}(R)$.  For more information on Chevalley groups, structure constants, and the various types of root systems $\Phi$, we refer the reader to \cite{VaPl96}. \par
Suppose $\Phi'$ is a root system generated by a subgraph of the Dynkin diagram of $\Phi$. The inclusion of root systems $\Phi'\subset \Phi$ gives rise to obvious functorial maps $\F(\Phi',R)\to\F(\Phi,R)$ for $\F=\G,\E,\St, \K_1$ or $\K_2$. In the classical case, where $\Phi=\Phi_r$ is of type $\A_r,\B_r,\C_r~\text{or}~\D_r$, one may take the direct limit of the above functor $\F$ over the system of maps induced by embedding $\Phi_r\subset \Phi_{r+1}$ to obtain the stable functor $$\F(R):=\F(\Phi_{\infty},R)=\lim_{\to r}\F(\Phi_r,R).$$ \par
In this paper, we will be mostly interested in the case where the Dynkin diagram of $\Phi$ is {\it simply-laced}, i.e., does not contain double bonds. This is equivalent to saying that all roots of $\Phi$ have the same length (i.e., are of type $\A,\D$ or $\E$). In this case, the defining relations of $\St(\Phi,R)$ take the following simpler form
\begin{align}
  &    x_\alpha(u) \, x_\alpha(v) = x_\alpha(u+v), \qquad u, v \in R,~ \alpha\in\Phi;\\[4pt]
   &    [x_\alpha(u), x_\beta(v)] = x_{\alpha + \beta}(C_{\alpha\beta}\, u v), \qquad u, v \in R,~\alpha+ \beta \in \Phi, \\[4pt]
  & [x_\alpha(u), x_\beta(v)] = 1, \qquad u, v \in R,~\alpha+ \beta \notin \Phi\cup 0.
\end{align}
\begin{rem}\label{rem_rel_ortk2}
    We remark that the functor $\K_{2}(\A_{r},R)$ coincides with the unstable linear functor $\K_{2,r+1}(R)$. On the other hand, the functor $\K_{2}(\D_{r},R):=\KSpin_{2,2r}(R)=\Ker(\St_{2}(\D_{r},R)\to\Spin_{2r}(R))$ differs from the unstable orthogonal functor $\KO_{2,2r}(R)=\Ker(\St_{2}(\D_{r},R)\to\Or_{2r}(R))$. However, they are related via the following exact sequence (see \cite[Section 2.3]{LaSi20} and \cite[Section 4.3]{Ba74}):
    \begin{center}
        \begin{tikzcd}
1 \arrow[r] & \K_{2}(\D_{r},R) \arrow[r]           & \KO_{2,2r}(R) \arrow[r]                 & \mu_{2}(R) \arrow[r]                 & 1,
\end{tikzcd}
    \end{center}
    where the group $\mu_{2}(R)=\{u\in R^{\times}~|~u^{2}=1\}$ is a central subgroup of $\E(\D_{r},R)$.
\end{rem}
For an ideal $I$ of $R$, the subgroup $\St(\Phi,I)$ is spanned by symbols $x_\alpha(s)$, where $\alpha \in \Phi$ and $s \in I$. We denote $\overline{\St}(\Phi,R,I)$ as the kernel  $\Ker(\St(\Phi,R)\twoheadrightarrow\St(\Phi,R/I))$. Obviously, $\overline{\St}(\Phi,R,I)$ coincides with the normal closure of $\St(\Phi,I)$ in $\St(\Phi,R)$.
\begin{defn} 
  An ideal $I$ is called a {\it splitting ideal} if the canonical projection ring homomorphism $\pi:R \twoheadrightarrow R/I$ splits, i.e., there exists $\sigma:R/I\to R$ such that $\pi\sigma=1$. It is noted that $XR[X]$ is a splitting ideal of $R[X]$.   
\end{defn}
From \cite[Section 3]{Si16} by definition, the {\it relative Steinberg group} $\St(\Phi,R,I)$ is the result of application of Loday’s relativization to the functor $\St(\Phi,-)$. If $I$ is a splitting ideal, then $\St(\Phi,R,I)$ coincides with $\overline{\St}(\Phi,R,I)$ (cf. \cite[Lemma 8]{Si16}). For the purposes of the present text, we shall focus primarily on the relative Steinberg groups corresponding to splitting ideals $XR[X]$ of $R[X]$. Further information on relative Steinberg groups can be found in \cite{LaSi20,Si16}.
 \begin{rem}\label{rem_split_k2}
     If $I$ a splitting ideal of $R$, then applying the functor $\K_{2}(\Phi,-)$ to $\pi$ and $\sigma$, we get \begin{center}
        \begin{tikzcd}
1 \arrow[r] & \K_{2}(\Phi,R,I) \arrow[r]           & \K_{2}(\Phi,R) \arrow[r, "\pi^{\Phi}", shift left]                 & \K_{2}(\Phi,R/I) \arrow[r]  \arrow[l, "\sigma^{\Phi}", shift left]               & 1.
\end{tikzcd}
    \end{center}
    Now, it is straightforward to verify that $\K_{2}(\Phi,R)$ decomposes as the semidirect product $$\K_{2}(\Phi,R/I)\ltimes \K_{2}(\Phi,R,I)$$ as a group. Consequently, there is a group homomorphism $$\varphi:\K_{2}(\Phi,R/I)\to\Aut(\K_{2}(\Phi,R,I))$$ given by conjugation via the splitting map $\sigma^{\Phi}$; explicitly, $\varphi(h)(g):=\sigma^{\Phi}(h)g\sigma^{\Phi}(h)^{-1}$.
 \end{rem}
\textbf{Notations:} For a ring $R$, we denote by $\lambda_{a}:R\to R_{a}$ the localization homomorphisms at an element $a\in R$. For a ring homomorphism $f: R\to S$, we denote by $f^{\Phi}$ the corresponding homomorphism of groups $\F(\Phi,R)\to\F(\Phi,S)$ for $\F=\G,\E,\St,\K_1$ or $\K_2$, provided the context is clear. If $\Phi'$ is an arbitrary irreducible simply-laced root system of rank $\geq 3$ embedded into another irreducible simply-laced root system $\Phi$, we let $l_{R}$ denote the corresponding homomorphism of Steinberg groups $\St(\Phi',R)\to\St(\Phi,R)$ and $\K_2$-functors $\K_{2}(\Phi',R)\to\K_{2}(\Phi,R)$ induced by this embedding.\par
\section{Monic Inversion Principles}\label{sec_mip}
In this section, we assume that $R$ is an integral domain and the root system $\Phi$ is of type $\A_{\geq 4},~\D_{\geq 5}$ or $\E_{6,7,8}$. The following result is a generalization of \cite[Proposition 4.3]{Tu83} and is essentially contained in Section $5$ of \cite{LaSi20, Si25}.
\begin{prop}\label{prop_inj_decomp}
    \begin{enumerate}
        \item The following homomorphisms are injective
    \begin{center}
    \begin{tikzcd}
\St(\Phi,R[X]) \arrow[r, "\lambda_{X}^{\Phi}"] & \St(\Phi,R[X,X^{-1}]) & \St(\Phi,R[X^{-1}]) \arrow[l, "\lambda_{X^{-1}}^{\Phi}"'].
\end{tikzcd}
\end{center} Moreover, the intersection of their images in $\St(\Phi,R[X,X^{-1}])$ coincides with the image of $\St(\Phi,R)$.
\item Let $\alpha\in\K_{2}(\Phi,R[X],XR[X])$. If $\lambda_{X}^{\Phi}(\alpha)$ in $\St(\Phi,R[X,X^{-1}])$ can be expressed as $l_{R[X,X^{-1}]}(\beta)\cdot\lambda^{\Phi}_{X^{-1}}(\gamma)$ for some $\beta\in\St(\Phi',R[X,X^{-1}])$ and $\gamma\in\St(\Phi,R[X^{-1}])$, then
$\alpha=l_{R[X]}(\alpha')$ for some $\alpha'\in\K_{2}(\Phi',R[X])$.
    \end{enumerate} 
\end{prop}
\begin{proof}
    $(1)$ This is a direct consequence of \cite[Theorem 3]{LaSi20} and \cite[Theorem 2]{LaSi17}. It should be noted that, although the hypothesis of \cite[Theorem 3]{LaSi20} requires the extra assumption that the canonical homomorphism $j_{R[X,X^{-1}]}$ (in the notation of that paper) be injective, however this condition was already verified later in the proof of \cite[Lemma 4.5]{Si25} (cf. also \cite[page 14]{Si25}). 
    % Let $\alpha\in\St(\Phi,R[X])$ such that $\lambda_{X}^{\Phi}(\alpha)=1$. Then for all maximal ideals $\mathfrak{m}$ of $R$, we have $\lambda_{X}^{\Phi}(\lambda_{\mathfrak{m}}^{\Phi}(\alpha))=\lambda_{\mathfrak{m}}^{\Phi}(\lambda_{X}^{\Phi}(\alpha))=1$. It follows from Proposition $\ref{prop_inj_decomp_local}(1)$ that $\lambda_{\mathfrak{m}}^{\Phi}(\alpha)=1$; hence, by Theorem $\ref{thm_lgp_st_slrs}(1)$, we get $\alpha=1$. Analogously, it can be shown that $\lambda_{X^{-1}}^{\Phi}$ is injective.\par
    % Let $\alpha\in\St(\Phi,R[X])$ and $\beta\in\St(\Phi,R[X^{-1}])$ such that $\lambda_{X}^{\Phi}(\alpha)=\lambda_{X^{-1}}^{\Phi}(\beta)$. If $\alpha'=\alpha\cdot\alpha^{-1}(0)$ and $\beta'=\beta\cdot\alpha^{-1}(0)$, then $\lambda_{X}^{\Phi}(\alpha')=\lambda_{X^{-1}}^{\Phi}(\beta')$. Note that $\alpha'\in\St(\Phi,R[X],XR[X])$, and hence, $\lambda_{\mathfrak{m}}^{\Phi}(\alpha')\in\St(\Phi,R_{\mathfrak{m}}[X],XR_{\mathfrak{m}}[X])$ for all maximal ideals $\mathfrak{m}$ of $R$. Since $$\lambda_{X}^{\Phi}(\lambda_{\mathfrak{m}}^{\Phi}(\alpha'))=\lambda_{X^{-1}}^{\Phi}(\lambda_{\mathfrak{m}}^{\Phi}(\beta')),$$ in view of Proposition $\ref{prop_inj_decomp_local}(1)$, the image $\lambda_{\mathfrak{m}}^{\Phi}(\alpha')$ is trivial. Then again by Theorem $\ref{thm_lgp_st_slrs}(1)$, we get $\alpha'=1$, and hence $\alpha=\alpha(0)\in\St(\Phi,R)$. Analogously, it can be shown that $\beta\in \St(\Phi,R)$.
    \par
    $(2)$ From \cite[Lemma 5.6]{Si25}, we have $l_{R[X]}(\alpha')=\alpha$ for some $\alpha'\in \St(\Phi',R[X])$. Then $\eta^{\Phi} l_{R[X]}(\alpha')=1$ in $\E(\Phi,R[X])$. Consider the following commutative diagram: 
\begin{center}
        \begin{tikzcd}
\St(\Phi,R[X]) \arrow[r,"\eta^{\Phi}"]                 & \E(\Phi,R[X])                 \\
\St(\Phi',R[X]) \arrow[r,"\eta^{\Phi'}"] \arrow[u, "l_{R[X]}"'] & \E(\Phi',R[X]) \arrow[u, "\theta"'],
\end{tikzcd}
    \end{center}
    where the right vertical map $\theta$ is injective (cf. \cite[page 80]{St78}). Due to the above diagram, we have $\theta\eta^{\Phi'}(\alpha')=1$ and hence $\eta^{\Phi'}(\alpha')=1$. Thus $\alpha'\in\K_{2}(\Phi',R[X])$. 
\end{proof}
% \begin{lem}
%     Let $\Phi$ be any simply-laced root system of rank $r \geq 3$ and $R$ be a Noetherian domain. Let $f\in R[X]$ be a monic polynomial. Let $\alpha\in\St(\Phi,R[X])$.
%     \begin{enumerate}
%         \item If $\lambda_{f}^{\Phi}(\alpha)=1$ in $\St(\Phi,R[X]_f)$, then $\alpha=1$.
%         \item If $\alpha\in\K_{2}(\Phi,R[X],XR[X])$ and $\lambda_{f}^{\Phi}(\alpha)$ lies in the image of $l_{R[X]_{f}}$, then $\alpha$ lies in the image of $l_{R[X]}$.
%     \end{enumerate}
% \end{lem}
% \begin{proof}
%     Write $f=X^n+\dots+a_n$ and $g=X^{-n}f$, then $(X^{-1},g)=1$. Denote $A=R[X,X^{-1}]$. Then $A_f=R[X^{-1}]_{X^{-1}g}$.
% \end{proof}
The following result is called the Zariski excision property for $\K_{2}$-functor. 
% \begin{prop}\label{prop_excision_st}
%      Let $R$ be a domain, and let $a,b\in R$ be a pair of coprime elements. Let $\Phi$ be a simply-laced type root system. Then the following assertions are true: 
%     \begin{itemize}
%         \item[(1)] If $\alpha\in\St(\Phi,R_{ab})$, then $\alpha$ can be written as $\lambda_{b}^{\Phi}(x)\cdot\lambda_{a}^{\Phi}(y)$ for some $x\in\St(\Phi,R_{a})$ and $y\in\St(\Phi,R_{b})$.
%         \item[(2)] If $x\in\St(\Phi,R_{a})$ and $y\in\St(\Phi,R_{b})$ such that $\lambda_{b}^{\Phi}(x)=\lambda_{a}^{\Phi}(y)$ in the group $\St(\Phi,R_{ab})$, then there exists $z\in\St(\Phi,R)$ such that $x=\lambda_{a}^{\Phi}(z)$ and $y=\lambda_{b}^{\Phi}(z)$. 
%     \end{itemize}
% \end{prop}
% \begin{proof}
%     See \cite[Lemma 5.7]{Si25} (more general Nisnevich gluing for $\K_{2}(\Phi,R)$ one can see \cite[Section 4]{LaSiVo23}).
% \end{proof}
\begin{prop}\label{prop_excision_k2}
Let $a,b\in R$ be a pair of coprime elements. Then the following assertions are true: 
\begin{enumerate}
\item[(1)] Assume that $\K_1(\Phi,R)=1$. For every $z\in \K_2(\Phi,R_{ab})$ there exist 
\(x\in \K_2(\Phi,R_a)\) and \(y\in \K_2(\Phi,R_b)\) such that
$z \;=\; \lambda_b^{\Phi}(x)\cdot\lambda_a^{\Phi}(y)$.
\item[(2)] Let \(x\in \K_2(\Phi,R_a)\) and \(y\in \K_2(\Phi,R_b)\) satisfy $\lambda_b^{\Phi}(x) = \lambda_a^{\Phi}(y)$ in $ \K_2(\Phi,R_{ab})$. Then there exists \(z\in \K_2(\Phi,R)\) such that $x = \lambda_a^{\Phi}(z)$ and $y = \lambda_b^{\Phi}(z)$.
\end{enumerate}
\end{prop}
\begin{proof}
$(1)$ Let $z\in \K_2(\Phi,R_{ab})\subseteq\St(\Phi,R_{ab})$. By \cite[Lemma 5.7]{Si25}, there exist $\tilde{x}\in\St(\Phi,A_a)$ and $ \tilde{y}\in\St(\Phi,A_b)$ such that
$z = \lambda_b^{\Phi}(\tilde{x})\cdot\lambda_a^{\Phi}(\tilde{y}) ~\text{in } \St(\Phi,R_{ab}).$ Since $\eta^{\Phi}(z)=1$, we have $
\lambda_b^{\Phi}(\eta^{\Phi}(\tilde{x}^{-1})) = \lambda_a^{\Phi}(\eta^{\Phi}(\tilde{y})) ~\text{in } \E(\Phi,R_{ab}).$\par
Note that the Chevalley group $\G(\Phi,R)$ is, in fact, a scheme (see \cite{Ch61}), and every scheme is a sheaf for the Zariski topology. This gives rise to the exact sequence \begin{center}
        \begin{tikzcd}
1 \arrow[r] & \G(\Phi,R) \arrow[r]           & \G(\Phi, R_{a})\times \G(\Phi, R_{b}) \arrow[r]                 & \G(\Phi, R_{ab}),
\end{tikzcd}
\end{center} where $\Spec(R_{ab})=\Spec(R_a)\cap \Spec(R_b)$ as a topological space. Thus, there exists a \(z'\in \G(\Phi,R)\) such that $\eta^{\Phi}(\tilde{x}^{-1}) = \lambda_a^{\Phi}(z'), \eta^{\Phi}(\tilde{y}) = \lambda_b^{\Phi}(z').$ Since $\K_1(\Phi,R)=1$, the homomorphism $\eta^{\Phi}:\St(\Phi,R)\to \G(\Phi,R)$ is surjective. Hence, we choose $\tilde{z}\in\St(\Phi,R)$ such that $\eta^{\Phi}(\tilde{z})=z'$. Then it is easy to verify that $\eta^{\Phi}(\tilde{x}\lambda_{a}^{\Phi}(\tilde{z}))=1$ and $\eta^{\Phi}(\lambda_{b}^{\Phi}(\tilde{z}^{-1})\tilde{y})=1$. Set $x := \tilde{x}\lambda_a^{\Phi}(\tilde{z}) \in \K_2(\Phi,R_a)$ and $y:=\lambda_b^{\Phi}(\tilde{z})^{-1}\tilde{y} \in \K_2(\Phi,R_b)$. Therefore, we get the required expression
\[
\begin{aligned}
z &= \lambda_b^{\Phi}(\tilde{x})\cdot\lambda_a^{\Phi}(\tilde{y}) \\
&= \lambda_b^{\Phi}(x\,\lambda_a^{\Phi}(\tilde{z})^{-1})\cdot
   \lambda_a^{\Phi}(\lambda_b^{\Phi}(\tilde{z})y) \\
&= \lambda_b^{\Phi}(x)\cdot\lambda_b^{\Phi}(\lambda_a^{\Phi}(\tilde{z}^{-1}))\cdot
   \lambda_a^{\Phi}(\lambda_b^{\Phi}(\tilde{z}))\cdot\lambda_a^{\Phi}(y)\\
&= \lambda_b^{\Phi}(x)\cdot\lambda_a^{\Phi}(y).   
\end{aligned}
\]\par
$(2)$ By \cite[Lemma 5.7]{Si25}, there exists 
\(u\in\St(\Phi,R)\) with $x = \lambda_a^{\Phi}(u), y = \lambda_b^{\Phi}(u)$. Then in $\E(\Phi,R_{a})$, we have $\lambda^{\Phi}_{a}(\eta^{\Phi}(u))=\eta^{\Phi}(\lambda^{\Phi}_{a}(u))=1$. Since the map \(\lambda_a^{\Phi}: \E(\Phi,R)\to \E(\Phi,R_a)\) is injective (cf. \cite[page 1280]{Ab83}) because 
\(R\hookrightarrow R_a\) is injective, we get $\eta^{\Phi}(u)=1$ in $\E(\Phi,R)$. This implies $u\in \K_2(\Phi,R)$. So, by  taking $z=u$, we get the assertion of $(2)$.
\end{proof}
\begin{thm}\label{thm_inj_surj_poly_monic}
     Let $f\in R[X]$ be a monic polynomial. 
     \begin{itemize}
         \item[(1)] If the homomorphism $\K_{2}(\Phi',R[X]_f)\to \K_{2}(\Phi,R[X]_f)$ is injective, then the homomorphism $\K_{2}(\Phi',R[X])\to \K_{2}(\Phi,R[X])$ is injective.
         \item[(2)] If homomorphisms $\K_{2}(\Phi',R)\to \K_{2}(\Phi,R)$ and $\K_{2}(\Phi',R[X]_f)\to \K_{2}(\Phi,R[X]_f)$ are surjective, then the homomorphism $\K_{2}(\Phi',R[X])\to \K_{2}(\Phi,R[X])$ is surjective. 
     \end{itemize}
\end{thm}
\begin{proof}
Write $f=X^n+a_{1}X^{n-1}+\dots+a_n$ and $g=X^{-n}f$, then $(X^{-1},g)=1$. Denote $A=R[X,X^{-1}]$. Then $A_f=R[X^{-1}]_{X^{-1}g}$. Consider the following commutative diagram: 
\begin{center}
 \begin{tikzcd}
\K_{2}(\Phi,R[X]) \arrow[d, hook] \arrow[r]     & \K_{2}(\Phi,R[X]_f) \arrow[d, hook] & \K_{2}(\Phi',R[X]_f) \arrow[d, hook] \arrow[r]  \arrow[l]    & \K_{2}(\Phi',R[X]) \arrow[d, hook]                \\
\St(\Phi,R[X]) \arrow[r, "\lambda^{\Phi}_{f}"] \arrow[d, "\lambda^{\Phi}_{X}"] & \St(\Phi,R[X]_f) \arrow[d,"\lambda^{\Phi}_{X}"]            & \St(\Phi',R[X]_f) \arrow[l, "l_{R[X]_f}"'] \arrow[d,"\lambda^{\Phi'}_{X}"] & \St(\Phi',R[X]) \arrow[l,"\lambda^{\Phi'}_{f}"'] \arrow[d,"\lambda^{\Phi'}_{X}"] \\
\St(\Phi,A) \arrow[r,"\lambda^{\Phi}_{f}"]                     & \St(\Phi,A_f)                           & \St(\Phi',A_{f}) \arrow[l,"l_{A_f}"']                      & \St(\Phi',A) \arrow[l,"\lambda^{\Phi'}_{f}"']                      \\
\St(\Phi,R[X^{-1}]) \arrow[u,"\lambda^{\Phi}_{X^{-1}}"'] \arrow[r]           & \St(\Phi,R[X^{-1}]_g) \arrow[u,"\lambda^{\Phi}_{X^{-1}}"']                 & \St(\Phi',R[X^{-1}]_g) \arrow[u," \lambda^{\Phi'}_{X^{-1}}"'] \arrow[l,"l_{R[X^{-1}]_g}"']            & \St(\Phi',R[X^{-1}]) \arrow[u,"\lambda^{\Phi'}_{X^{-1}}"'] \arrow[l, "\lambda^{\Phi'}_{g}"']   
% \arrow[phantom, from=1-1, to=2-2, "\scriptstyle I" description]
% \arrow[phantom, from=1-2, to=2-3, "\scriptstyle II" description]
% \arrow[phantom, from=1-3, to=2-4, "\scriptstyle III" description]
% \arrow[phantom, from=2-1, to=3-2, "\scriptstyle IV" description]
% \arrow[phantom, from=2-2, to=3-3, "\scriptstyle V" description]
% \arrow[phantom, from=2-3, to=3-4, "\scriptstyle VI" description]
\arrow[phantom, from=3-1, to=4-2, "\scriptstyle \text{I}" description]
%\arrow[phantom, from=3-2, to=4-3, "\scriptstyle VIII" description]
\arrow[phantom, from=3-3, to=4-4, "\scriptstyle \text{II}" description].
\end{tikzcd}
\end{center}
   $(1)$ Let $\alpha\in \K_{2}(\Phi',R[X])$ such that $l_{R[X]}(\alpha)=1$ in $\K_{2}(\Phi,R[X])$. Then $l_{R[X]_f}\lambda_{f}^{\Phi'}(\alpha)=\lambda_{f}^{\Phi}l_{R[X]}(\alpha)=1$. By hypothesis, we have $\lambda_{f}^{\Phi'}(\alpha)=1$, then $\lambda_{X}^{\Phi'}\lambda_{f}^{\Phi'}(\alpha)=1$. This can be written as $\lambda_{f}^{\Phi'}(\lambda_{X}^{\Phi'}(\alpha))=\lambda_{X^{-1}}^{\Phi'}(1)$. Then apply Proposition $\ref{prop_excision_k2}(2)$ on (II), we get $\beta\in \St(\Phi',R[X^{-1}])$ such that $\lambda_{X^{-1}}^{\Phi'}(\beta)=\lambda_{X}^{\Phi'}(\alpha)$ and $\lambda_{g}^{\Phi'}(\beta)=1$. It follows from Proposition $\ref{prop_inj_decomp}(1)$ that $\beta\in\St(\Phi',R)$. Since the composite of the homomorphisms 
    \begin{center}
        \begin{tikzcd}
 \St(\Phi',R)\arrow[r] & \St(\Phi',R[X^{-1}]_g) \arrow[r, "\tau"] & \St(\Phi',R)
\end{tikzcd}
    \end{center} coincides with the identity map (here $\tau$ is a map evaluated at $0$), we conclude that $\beta=\tau(\lambda_{g}^{\Phi'}(\beta))=1$. Then $\lambda_{X}^{\Phi'}(\alpha)=1$. Since $\lambda_{X}^{\Phi'}$ is injective, we have $\alpha=1$. This proves the injectivity.\par
$(2)$ Let $\alpha\in \K_{2}(\Phi,R[X])$. Since $XR[X]$ is a splitting ideal of $R[X]$, following Remark $\ref{rem_split_k2}$, we can decompose $\alpha$ as $\alpha=\alpha_{1}\cdot\alpha_{2}$, where $\alpha_{1}\in\K_{2}(\Phi,R[X],XR[X])$ and $\alpha_{2}\in\K_{2}(\Phi,R)$. Since $\lambda^{\Phi}_{f}(\alpha)$ belongs to the image of $l_{R[X]_{f}}$, it follows that $\lambda^{\Phi}_{f}(\alpha_{1})$ belongs to the image of $l_{R[X]_{f}}$. Then there exists $\delta\in\K_{2}(\Phi',R[X]_f)$ such that $\lambda^{\Phi}_{f}(\alpha_{1})=l_{R[X]_f}(\delta)$. Now, apply \cite[Lemma 5.7(1)]{Si25} on (II), we get $$\lambda^{\Phi'}_{X}(\delta)=\lambda^{\Phi'}_{f}(\beta)\cdot \lambda^{\Phi'}_{X^{-1}}(\gamma)$$ for some $\beta\in\St(\Phi',A)$ and $\gamma\in\St(\Phi',R[X^{-1}]_g)$. From the commutativity of the diagram, it is easy to verify the following: \[l_{A_{f}}\lambda_{f}^{\Phi'}=\lambda_{f}^{\Phi}l_{A},~ l_{A_{f}}\lambda_{X}^{\Phi'}=\lambda_{X}^{\Phi}l_{R[X]_{f}},~l_{A_{f}}\lambda_{X^{-1}}^{\Phi'}=\lambda_{X^{-1}}^{\Phi}l_{R[X^{-1}]_{g}} ~\text{and}~ \lambda_{X}^{\Phi}\lambda_{f}^{\Phi}=\lambda_{f}^{\Phi}\lambda_{X}^{\Phi}.\] Then, inside the image of $l_{*}$, we have $$\lambda^{\Phi}_{f}(l_{A}(\beta)^{-1}\cdot\lambda_{X}^{\Phi}(\alpha_{1}))=\lambda^{\Phi}_{X^{-1}}(l_{R[X^{-1}]_g}(\gamma)).$$ Apply \cite[Lemma 5.7(2)]{Si25} on (I), there exists $\omega\in \St(\Phi,R[X^{-1}])$ such that $\lambda_{X}^{\Phi}(\alpha_{1})=l_{A}(\beta)\cdot\lambda^{\Phi}_{X^{-1}}(\omega)$. Then by Proposition $\ref{prop_inj_decomp}(2)$, for some $\alpha_{1}'\in\K_{2}(\Phi',R[X])$, we have $\alpha_{1}=l_{R[X]}(\alpha_{1}')$. Also, by hypothesis, for some $\alpha_{2}'\in\K_{2}(\Phi',R)$, we have $\alpha_{2}=l_{R}(\alpha_{2}')$. Thus, the element $\alpha_{1}'\cdot\alpha_{2}'$ is a preimage of $\alpha$. This proves the surjectivity.
\end{proof}
Recall that a monic polynomial $f(X)\in R[X]$ is called a special monic if $f(0)=1$. The following is an analog of monic inversion principles in Laurent polynomial rings.
\begin{thm}\label{thm_inj_surj_Laur_poly_spmonic}
    Let $f\in R[X]$ be a special monic polynomial and  $A=R[X,X^{-1}]$.
    \begin{enumerate}
        \item If the homomorphism $\K_{2}(\Phi',A_f)\to \K_{2}(\Phi,A_f)$ is injective, then the homomorphism $\K_{2}(\Phi',A)\to \K_{2}(\Phi,A)$ is injective.
        \item Assume that $\K_{1}(\Phi,R[X])$ is trivial. If homomorphisms $\K_{2}(\Phi',R)\to \K_{2}(\Phi,R)$ and $\K_{2}(\Phi',A_f)\to \K_{2}(\Phi,A_f)$ are surjective, then the homomorphism $\K_{2}(\Phi',A)\to \K_{2}(\Phi,A)$ is surjective.
    \end{enumerate}
\end{thm}
\begin{proof}
Since $f$ is a special monic polynomial, $X$ and $f$ are comaximal elements of $R[X]$. We consider the following commutative diagram:
    \begin{center}
    \begin{tikzcd}
%\K_{2}(\Phi,R[X]) \arrow[r] \arrow[d, hook] & \K_{2}(\Phi,R[X]_{f}) \arrow[d, hook] & \K_{2}(\Phi',R[X]_{f}) \arrow[d, hook] \arrow[l] & \K_{2}(\Phi',R[X]) \arrow[d, hook] \arrow[l] \\
\K_{2}(\Phi,R[X]) \arrow[r,"\lambda^{\Phi}_{f}"] \arrow[d,"\lambda^{\Phi}_{X}"]                & \K_{2}(\Phi,R[X]_{f})        \arrow[d,"\lambda^{\Phi}_{X}"]         & \K_{2}(\Phi',R[X]_{f}) \arrow[l,"l_{R[X]_f}"']  \arrow[d,"\lambda^{\Phi'}_{X}"]              & \K_{2}(\Phi',R[X]) \arrow[l,"\lambda^{\Phi'}_{f}"']  \arrow[d,"\lambda^{\Phi'}_{X}"]               \\
\K_{2}(\Phi,A) \arrow[r,"\lambda^{\Phi}_{f}"]       & \K_{2}(\Phi,A_{f})        & \K_{2}(\Phi',A_{f})  \arrow[l,"l_{A_f}"']       & \K_{2}(\Phi',A)  \arrow[l,"\lambda^{\Phi'}_{f}"'] .
%\\ \K_{2}(\Phi,A) \arrow[u, hook] \arrow[r] & \K_{2}(A_{f}) \arrow[u, hook] & \K_{2}(\Phi',A_{f}) \arrow[u, hook] \arrow[l] & \K_{2}(\Phi',A) \arrow[u, hook] \arrow[l]
 \arrow[phantom, from=1-1, to=2-2, "\scriptstyle \text{I}" description]
%\arrow[phantom, from=1-2, to=2-3, "\scriptstyle II" description]
\arrow[phantom, from=1-3, to=2-4, "\scriptstyle \text{II}" description]
% \arrow[phantom, from=2-1, to=3-2, "\scriptstyle I" description]
% \arrow[phantom, from=2-2, to=3-3, "\scriptstyle V" description]
% \arrow[phantom, from=2-3, to=3-4, "\scriptstyle II" description]
% \arrow[phantom, from=3-1, to=4-2, "\scriptstyle VII" description]
%  \arrow[phantom, from=3-2, to=4-3, "\scriptstyle VIII" description]
%  \arrow[phantom, from=3-3, to=4-4, "\scriptstyle IX" description]
\end{tikzcd}
\end{center}
$(1)$ Let $\alpha\in \K_{2}(\Phi',A)$ such that $l_{A}(\alpha)=1$ in $\K_{2}(\Phi,A)$. Then $l_{A_f}\lambda_{f}^{\Phi'}(\alpha)=\lambda_{f}^{\Phi}l_{A}(\alpha)=1$. By hypothesis, we have $\lambda_{f}^{\Phi'}(\alpha)=1=\lambda_{X}^{\Phi'}(1)$. Applying Proposition $\ref{prop_excision_k2}(2)$ on (II), we get $\beta\in\St_{r}(\Phi',R[X])$ such that $\lambda_{X}^{\Phi'}(\beta)=\alpha$ and $\lambda_{f}^{\Phi'}(\beta)=1$. Note that $Xf$ is a monic polynomial in $R[X]$ and $\lambda_{Xf}^{\Phi'}(\beta)=\lambda_{X}^{\Phi'}(\lambda_{f}^{\Phi'}(\beta))=1$. As in the last paragraph of proof of Theorem $\ref{thm_inj_surj_poly_monic}(1)$, we conclude $\beta=1$ and hence $\alpha=1$. This proves the injectivity.\par
$(2)$ Let $\alpha\in \K_{2}(\Phi,A)$. Since $\lambda_{f}^{\Phi}(\alpha)\in\K_{2}(\Phi,A_f)$, by hypothesis, there exists $\omega\in\K_{2}(\Phi',A_f)$ such that $\lambda_{f}^{\Phi}(\alpha)=l_{A_{f}}(\omega)$. Applying Proposition $\ref{prop_excision_k2}(2)$ on (II), we get $\omega=\lambda_{f}^{\Phi'}(\gamma)\cdot\lambda_{X}^{\Phi'}(\delta)$ for some $\gamma\in\K_{2}(\Phi',A)$ and $\delta\in\K_{2}(\Phi',R[X]_{f})$. Consequently, after applying $l_{*}$ and due to the commutativity, we have that $$\lambda_{f}^{\Phi}(\alpha\cdot l_{A}(\gamma)^{-1})=\lambda_{X}^{\Phi}(l_{R[X]_f}(\delta)).$$ Now, applying Proposition $\ref{prop_excision_k2}(1)$ on (I), we get $\beta\in\K_{2}(\Phi,R[X])$ such that $\alpha=\lambda_{X}^{\Phi}(\beta)\cdot l_{A}(\gamma)$ and $\lambda_{f}^{\Phi}(\beta)=l_{R[X]_f}(\delta)$. Since $Xf\in R[X]$ is a monic polynomial and $\lambda^{\Phi}_{Xf}(\beta)$ lies in the image of $l_{R[X]_{Xf}}$. The verbatim copy of the proof of Theorem $\ref{thm_inj_surj_poly_monic}(2)$ implies that $\beta=l_{R[X]}(\beta')$ for some $\beta'\in\K_{2}(\Phi',R[X])$. Hence $\lambda_{X}^{\Phi}(\beta)=\lambda_{X}^{\Phi}l_{R[X]}(\beta')=l_{A}(\beta'')$, where $\beta''=\lambda_{X^{-1}}^{\Phi'}(\beta')\in\K_{2}(\Phi',A)$. Therefore, $\alpha=l_{A}(\beta''\cdot\gamma)$. This proves the surjectivity.
\end{proof}

\section{Stability Theorems over Polynomial Rings}\label{sec_pr}
We begin with the following concept, introduced by H. Bass \cite[Chapter 5, Section 3]{Ba68}.
\begin{defn}
$A$ ring $R$ satisfies condition $(\SR_{m})$ if for every $(a_1, \ldots, a_m) \in \Um_m(R)$, there exists $x_1, \ldots, x_{m-1}\in R$ such that $(a_1 + a_m x_1, \cdots ,a_{m-1} + a_{m}x_{m-1})\in\Um_{m-1}(R)$. The condition $(\SR_m)$ implies $(\SR_{m+1})$ for every $m>0$. Moreover, for any $n\geq m$, the condition $(\SR_{m})$ implies $(\SR_{n})$.
\end{defn}
 The stable rank for a ring $R$ is the least $m$ such that $(\SR_{m+1})$ holds, and it is denoted by $\sr(R)$. Moreover, if $\sr(R)\leq k$, then $R$ satisfies the condition $(\SR_{k+1})$.\par
 %If $R$ is a Noetherian ring, then by \cite[Theorem 3.5]{Ba68} $\sr(R)\leq\dim(R)+1$.
  The following result is due to Bass \cite[Theorem 3.5]{Ba68}.
  \begin{thm}\label{thm_sr_bass}
      Let $R$ be a Noetherian ring. Then $\sr(R)\leq\dim(R)+1$.  
  \end{thm}
  The following stability result is due to Stein \cite[Corollary 3.2, Theorem 4.1]{St78} utilizing Theorem $\ref{thm_sr_bass}$.
\begin{thm}\label{thm_sr_stein}
      Let $R$ be a ring. Then the homomorphism $\K_{2}(\Phi_r,R)\to \K_{2}(\Phi_{r+1},R)$ is surjective for $r\geq \sr(R)$ if $\Phi_r=\A_{r\geq 1}$, $r\geq \sr(R)+1$ if $\Phi_r=\D_{r\geq 2}$, and $\sr(R)\leq r-2$ if $\Phi_{r}=\E_{r=5,6,7}$.
\end{thm}
The following are Dynkin diagrams for the simply-laced root systems $\A_{r\geq 4},~\D_{r\geq 4}$ and $\E_{r=6,7,8}$, each with $r$ vertices (note that $\E_4=\A_4$ and $\E_5=\D_5$).
\begin{figure}[H]
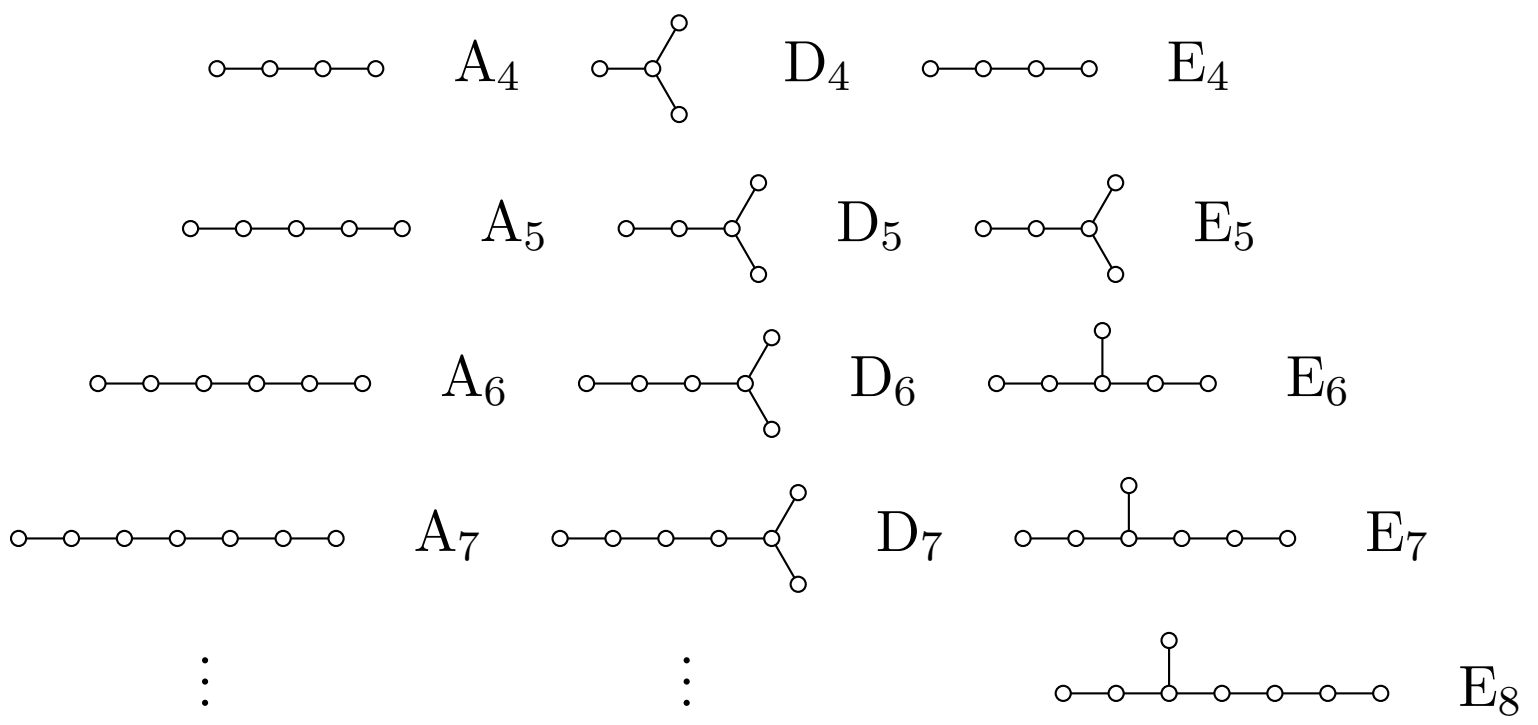

    \centering
\[
\dynkin[mark=o]{A}{4} \quad \text{A}_4 \quad \dynkin[mark=o]{D}{4} \quad \text{D}_4 \quad \dynkin[mark=o]{A}{4} \quad \text{E}_4
\]
\[
\dynkin[mark=o]{A}{5} \quad \text{A}_5 \quad \dynkin[mark=o]{D}{5} \quad \text{D}_5 \quad 
\dynkin[mark=o]{D}{5} \quad \text{E}_5
\]
\[
\dynkin[mark=o]{A}{6} \quad \text{A}_6 \quad \dynkin[mark=o]{D}{6} \quad \text{D}_6 \quad
\dynkin[mark=o]{E}{6} \quad \text{E}_6
\]
\[\dynkin[mark=o]{A}{7} \quad \text{A}_7 \quad \dynkin[mark=o]{D}{7} \quad \text{D}_7\quad \dynkin[mark=o]{E}{7} \quad \text{E}_7
\]
\[\quad \quad\quad\quad \quad\vdots \qquad\qquad\qquad\qquad  \vdots \quad\quad \qquad\quad \quad \dynkin[mark=o]{E}{8} \quad \text{E}_8
\]
    \caption{Dynkin Diagrams}
    \label{dynkin_diagrams}
\end{figure} 
In the following diagram, arrows indicate containments for the subsystem embedding of a simply-laced root system. The embeddings in ``dasharrow" come from subgraphs of the above Dynkin diagram. The (full rank) embeddings in ``squigarrow" are not visible as subdiagrams of the corresponding Dynkin diagrams; their existence follows from Dynkin's classification of maximal closed subsystems of maximal rank root systems obtained by removing one extreme fundamental root (see Borel–de Siebenthal classification \cite{BoSi49}). Consequently, all of these are root subsystem embeddings and therefore induce the corresponding homomorphisms of Steinberg groups.
\begin{figure}[H]
    \centering
    \begin{tikzcd}
\D_{4} \arrow[r,dashrightarrow] \arrow[rd, dashrightarrow]                                            & \D_{5} \arrow[r, dashrightarrow] \arrow[rd, dashrightarrow]                     & \D_{6} \arrow[r, dashrightarrow] \arrow[rd, dashrightarrow]             & \D_{7} \arrow[r, dashrightarrow] \arrow[rd, dashrightarrow]             & \D_{8} \arrow[r, dashrightarrow] \arrow[d, rightsquigarrow] & \cdots \\
\E_{4} \arrow[r, dashrightarrow]                                                       & \E_{5} \arrow[r, dashrightarrow] \arrow[u, equals] & \E_{6} \arrow[r, dashrightarrow]                        & \E_{7} \arrow[r, dashrightarrow]                        & \E_{8}           &        \\
\A_{4} \arrow[r, dashrightarrow] \arrow[ru, dashrightarrow] \arrow[ruu, dashrightarrow] \arrow[u, equals] & \A_{5} \arrow[r, dashrightarrow] \arrow[ru, dashrightarrow] \arrow[ruu, dashrightarrow]         & \A_{6} \arrow[r, dashrightarrow] \arrow[ru, dashrightarrow] \arrow[ruu, dashrightarrow] & \A_{7} \arrow[r, dashrightarrow] \arrow[u, rightsquigarrow]\arrow[ru, dashrightarrow] \arrow[ruu, dashrightarrow] & \A_{8} \arrow[r, dashrightarrow] \arrow[u, rightsquigarrow] & \cdots
\end{tikzcd}
    \caption{Root Embeddings}
    \label{root_embedding}
\end{figure}
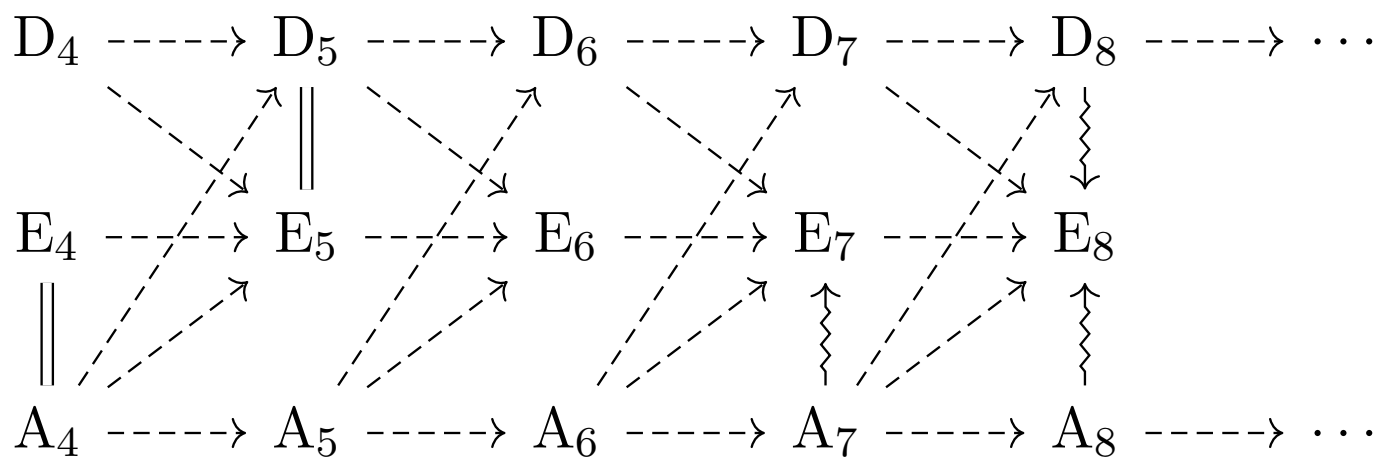
% \begin{center}
% \begin{tikzcd}
% \D_{4} \arrow[r,dashrightarrow] \arrow[rd, dashrightarrow]                                            & \D_{5} \arrow[r, dashrightarrow] \arrow[rd, dashrightarrow]                     & \D_{6} \arrow[r, dashrightarrow] \arrow[rd, dashrightarrow]             & \D_{7} \arrow[r, dashrightarrow] \arrow[rd, dashrightarrow]             & \D_{8} \arrow[r, dashrightarrow] \arrow[d, rightsquigarrow] & \cdots \\
% \E_{4} \arrow[r, dashrightarrow]                                                       & \E_{5} \arrow[r, dashrightarrow] \arrow[u, equals] & \E_{6} \arrow[r, dashrightarrow]                        & \E_{7} \arrow[r, dashrightarrow]                        & \E_{8}           &        \\
% \A_{4} \arrow[r, dashrightarrow] \arrow[ru, dashrightarrow] \arrow[ruu, dashrightarrow] \arrow[u, equals] & \A_{5} \arrow[r, dashrightarrow] \arrow[ru, dashrightarrow] \arrow[ruu, dashrightarrow]         & \A_{6} \arrow[r, dashrightarrow] \arrow[ru, dashrightarrow] \arrow[ruu, dashrightarrow] & \A_{7} \arrow[r, dashrightarrow] \arrow[u, rightsquigarrow]\arrow[ru, dashrightarrow] \arrow[ruu, dashrightarrow] & \A_{8} \arrow[r, dashrightarrow] \arrow[u, rightsquigarrow] & \cdots
% \end{tikzcd}
% \end{center}
It is noted that by deleting suitable roots, we obtain more closed subsystems (e.g., $\D_6+\A_1\subset\E_7,~\A_5+\A_2\subset\E_7,~\E_6+\A_2\subset\E_8,~\A_4+\A_4\subset\E_8$). These are reducible; since we are only interested in irreducible ones, we omit them. \par
The following results generalize \cite[Theorem 5.3]{Tu83}.
\begin{thm}\label{thm_surj_stab_poly}
     Let $R$ be a Noetherian domain and $A=R[X_{1},\ldots,X_{n}]$. Let $\Phi'\subset \Phi$ be one of the following embeddings of root systems:
    \begin{enumerate}
        \item[$(i)$] $\A_{r}\subset \A_{r+1}$, with $r\geq\max\{3,\dim(R)+1\}$;
        \item[$(ii)$] $\D_{r}\subset \D_{r+1}$, with $r\geq\max\{4,\dim(R)+2\}$;
        \item[$(iii)$] $\E_{r}\subset \E_{r+1}~(r=5,6,7)$, with $\dim(R)\leq r-3$;
        \item[$(iv)$] $\D_{r}\subset \E_{r+1}~(r=6,7)$ and $\D_{8}\subset \E_{8}$, with $\dim(R)\leq 2$;
        \item[$(v)$] $\A_{r}\subset \E_{r+1}~(r=5,6,7)$ and $\A_{r}\subset \D_{r+1}~(r\geq 5)$, with $\dim(R)\leq 1$;
         \item[$(vi)$] $\A_{r}\subset \E_{r}~(r=7,8)$, with $\dim(R)\leq 1$;
         \item[$(vii)$] $\A_{4}\subset\A_{\geq 4},~\A_{4}\subset\D_{\geq 5}$ and $\A_{4}\subset\E_{6,7,8}$, with $\dim(R)\leq 1$.
    \end{enumerate}
    Then the canonical homomorphism $\Gamma_{A}:\K_{2}(\Phi',A)\to \K_{2}(\Phi,A)$ is surjective under the given conditions. 
\end{thm}
\begin{proof}
    For (vii) see \cite[Theorem 5.9]{Si25}. We will prove the theorem by induction on $n$.\par
    \textbf{Step 1.} Assume $n=0$. By Theorem $\ref{thm_sr_bass}$, we have $\sr(R)\leq\dim(R)+1$. For (i)--(iii), our assumption on $\dim(R)$ guarantees that the map $\Gamma_{R}:\K_{2}(\Phi',R)\to \K_{2}(\Phi,R)$ is surjective by Theorem $\ref{thm_sr_stein}$.\par If $\dim(R)\leq 2$, by (iii) the composite maps of the following diagram are surjective (recall $\D_5=\E_5$).
    
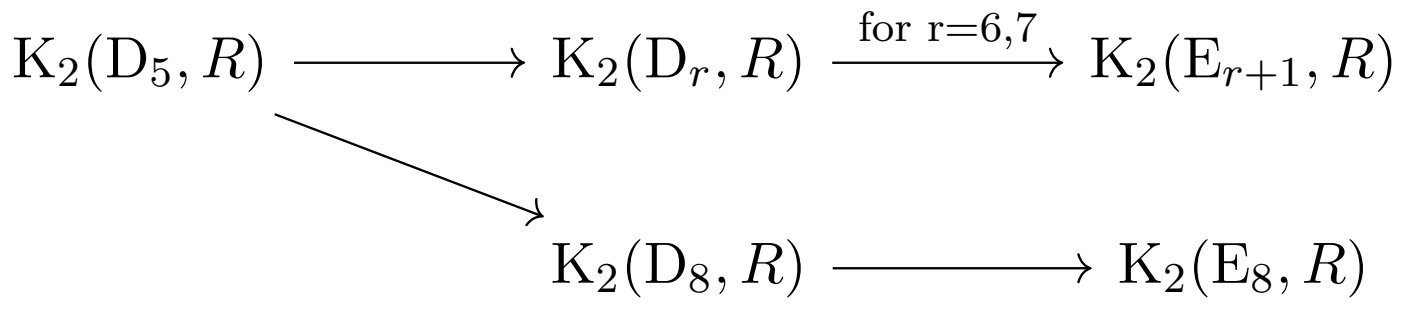
\begin{center}
        \begin{tikzcd}[column sep=4em]
\K_{2}(\D_5,R) \arrow[r]\arrow[rd]  &     \K_{2}(\D_{r},R) \arrow[r, "\text{for~r=6,7}"]         & \K_{2}(\E_{r+1},R)\\
      & \K_{2}(\D_{8},R)  \arrow[r] & \K_{2}(\E_{8},R)
\end{tikzcd}
    \end{center}
%     \begin{center}
%         \begin{tikzcd}[column sep=4em]
%                         & \K_{2}(\D_{r},R) \arrow[r, "\text{for~r=6,7}"] & \K_{2}(\E_{r+1},R) \\
% \K_{2}(\D_5,R) \arrow[ru] \arrow[rd] &                  &   \\
%                         & \K_{2}(\D_{8},R) \arrow[r]      & \K_{2}(\E_{8},R).
% \end{tikzcd}
%     \end{center}
Hence, the right maps are surjective, which shows that $\Gamma_{R}$ is surjective for (iv). For (v) and (vi), since $\dim(R)\leq 1$, by (vii), the composite maps of the following diagram are surjective.
%     \begin{center}
% \begin{tikzcd}
%             &                         & \K_{2}(\E_{r+1},R) \\
% \K_{2}(\A_4,R) \arrow[r] & \K_{2}(\A_{r},R) \arrow[ru, "\text{for~r=5,6,7}"'] \arrow[rd, "\text{for ~r$\geq$5}"] &   \\
%             &                         & \K_{2}(\D_{r+1},R).
% \end{tikzcd}
%     \end{center} 
    
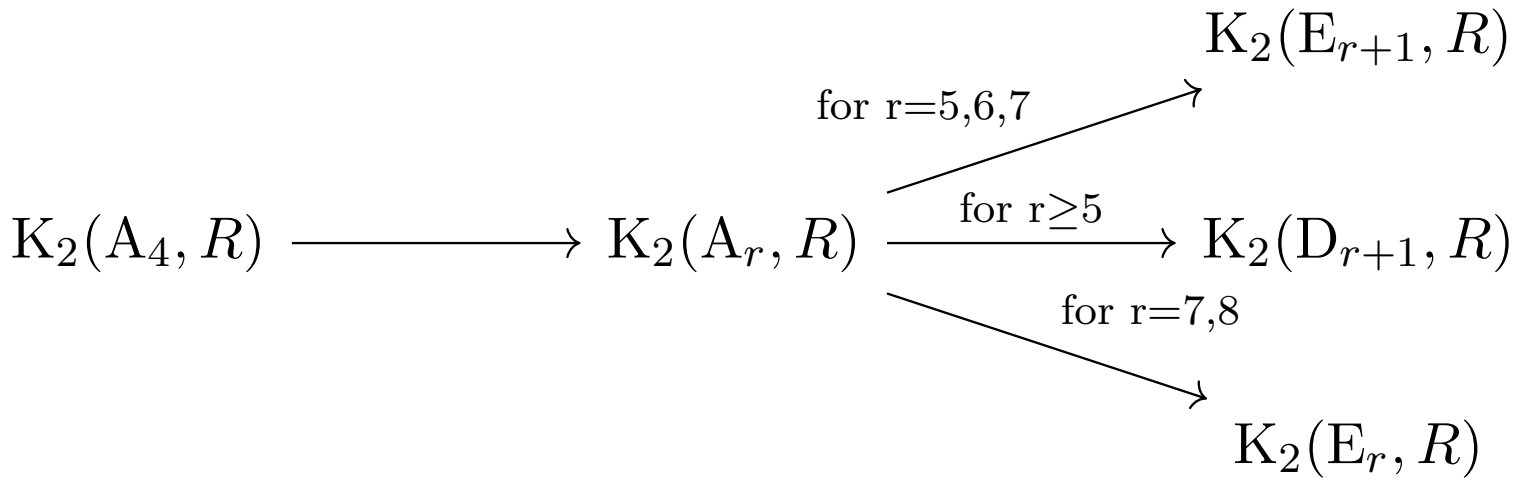
\begin{center}
        \begin{tikzcd}[column sep=5em]
            &                                                   & \K_{2}(\E_{r+1},R) \\
\K_{2}(\A_4,R) \arrow[r] & \K_{2}(\A_{r},R) \arrow[ru, "\text{for~r=5,6,7}"] \arrow[r, "\text{for ~r$\geq$5}"] \arrow[rd, "\text{for ~r=7,8}"] & \K_{2}(\D_{r+1},R) \\
            &                                                   & \K_{2}(\E_{r},R)
\end{tikzcd}
    \end{center}
Hence, the right maps are surjective, which shows that $\Gamma_{R}$ is surjective for (v) and (vi).
% For (vii); since $\dim(R)\leq 1$, by (v) the composite maps of the following diagram is surjective 
% \begin{center}
%   \begin{tikzcd}
% \K_{2}(\A_{6},R) \arrow[r]  &     \K_{2}(\A_{r},R) \arrow[r]         & \K_{2}(\E_{r},R)~~(\text{for}~r=7,8).
% \end{tikzcd}
% \end{center}
% Hence, the right one is surjective, which shows that $\Gamma_{R}$ is surjective for (vii).
% Similarly, for (viii), the map $\Gamma_{R}$ is surjective. Indeed, since $\dim(R)\leq 2$, by (iv) the composite maps of the following diagram is surjective
% \begin{center}
%   \begin{tikzcd}
% \K_{2}(\D_{7},R) \arrow[r]  &     \K_{2}(\D_{8},R) \arrow[r]         & \K_{2}(\E_{8},R).
% \end{tikzcd}
% \end{center}
% Hence, the right map is surjective. 
This yields the base case of the induction.\par   
\textbf{Step 2.} Assume that $n\geq 1$. If $B=R[X_{2},\ldots,X_{n}]$, then $A=B[X_{1}]$. Let $S$ be the multiplicative system in $A$ consisting of polynomials $f$ such that for sufficiently large $m$, the polynomial $f$ becomes monic with respect to $Y_{1}$ after substitutions $X_{1}:=Y_{1},X_{i}:=Y_{i}+Y_{1}^{m^{i-1}}$ for $i=2,\ldots,n$. From \cite[Lemma 6.2]{Su77b}, it is known that $\dim(S^{-1}A)\leq\dim(R)$ and hence, by Theorem $\ref{thm_sr_bass}$, $\sr(S^{-1}A)\leq\dim(R)+1$.\par For (i)--(iii), our assumption on $\dim(R)$ guarantees that $S^{-1}A$ satisfies the stable rank condition given in the hypothesis of Theorem $\ref{thm_sr_stein}$. Thus the map $\Gamma_{S^{-1}A}:\K_{2}(\Phi',S^{-1}A)\to \K_{2}(\Phi,S^{-1}A)$ is surjective. With the assumption on $\dim(R)$, it is easy to show that the map $\Gamma_{S^{-1}A}$ is surjective for (iv)--(vi) as in Step 1.\par
By \cite[Lemma 3.3]{LaSiVo23}, the functor $\K_{2}$ commutes with filtered colimits, then there exists $f\in S$ such that $\Gamma_{A_{f}}:\K_{2}(\Phi',A_f)\to \K_{2}(\Phi,A_f)$ is surjective. Also, by induction assumption, the map $\Gamma_{B}:\K_{2}(\Phi',B)\to \K_{2}(\Phi,B)$ is surjective. Hence, since $A_f=B[X_1]_f$, the assertion now follows from Theorem $\ref{thm_inj_surj_poly_monic}(2)$.
\end{proof}
The following results are due to Suslin--Tulenbaev \cite[Theorem 4.1]{SuTu76} and Panin \cite[Theorem 9.4]{Pa90} respectively.
\begin{thm}\label{thm_sutu_inj}
     Let $R$ be a ring. Then the homomorphism $\K_{2}(\A_r,R)\to \K_{2}(\A_{r+1},R)$ is injective for $r\geq \sr(R)+1$.
\end{thm}
\begin{thm}\label{thm_stb_dr_panin}
    Let $R$ be either a field, a principal ideal domain, or a Dedekind domain. Then the canonical homomorphism $\KO_{2,2r}(R)\to \KO_{2,2r+2}(R)$ is surjective for $r\geq 4$ and injective for $r\geq 5$.
\end{thm}
\begin{thm}\label{thm_inj_stab_poly}
Let $R$ be a Noetherian domain and $A=R[X_1,\ldots,X_n]$. Let $\Phi'\subset \Phi$ be one of the following embeddings of root systems:
    \begin{enumerate}
        \item[$(i)$] $\A_{r}\subset \A_{r+1}$ for $r\geq\max\{4,\dim(R)+2\}$;
        \item[$(ii)$] $\D_{r}\subset \D_{r+1}$ for $r\geq 5$.
    \end{enumerate}
In $(ii)$, assume that $R$ is either a principal ideal domain or a Dedekind domain containing a field. Then the canonical homomorphism $\Gamma_{A}:\K_{2}(\Phi',A)\to \K_{2}(\Phi,A)$ is injective under the given conditions. 
\end{thm}
\begin{proof}
    (i) Let $S$ be the multiplicative system in $A$ consisting of polynomials $f$ such that for sufficiently large $m$, the polynomial $f$ becomes monic with respect to $Y_{1}$ after substitutions $X_{1}:=Y_{1},~X_{i}:=Y_{i}+Y_{1}^{m^{i-1}}$ for $i=2,\ldots,n$. From \cite[Lemma 6.2]{Su77b}, it is known that $\dim(S^{-1}A)\leq\dim(R)$, and hence by Theorem $\ref{thm_sr_bass}$, $\sr(S^{-1}A)\leq\dim(R)+1$. By Theorem $\ref{thm_sutu_inj}$, the homomorphism $\Gamma_{A_f}:\K_{2}(A_r,A_f)\to \K_{2}(\A_{r+1},A_f)$ is injective for some monic polynomial $f$. The assertion now follows from Theorem $\ref{thm_inj_surj_poly_monic}$(1).\par
    (ii) Since a principal ideal domain and a Dedekind domain are regular rings, by \cite[Theorem 1.1]{Si25}, we have $\K_{2}(\D_r,R)\cong\K_{2}(\D_r,R[X_1,\ldots,X_n])$ for $r\geq 5$. Consider the following commutative diagram:
    \begin{center}
        \begin{tikzcd}
\K_{2}(\D_r,A) \arrow[r, "\Gamma_{A}"]                 & \K_{2}(\D_{r+1},A)                 \\
\K_{2}(\D_r,R) \arrow[r,"\Gamma_{R}"] \arrow[u, "\cong"'] & \K_{2}(\D_{r+1},R). \arrow[u, "\cong"']
\end{tikzcd}
    \end{center}
 To get the result, it is enough to prove that $\Gamma_{R}$ is injective. From Remark $\ref{rem_rel_ortk2}$, we get the following diagram of exact sequences:
 \begin{center}
        \begin{tikzcd}
1 \arrow[r] & \K_{2}(\D_{r+1},R) \arrow[r,"i"]           & \KO_{2,2r+2}(R) \arrow[r]                 & \mu_{2}(R) \arrow[r]                 & 1 \\
1 \arrow[r] & \K_{2}(\D_{r},R) \arrow[r,"j"] \arrow[u,"\Gamma_{R}"'] & \KO_{2,2r}(R) \arrow[r] \arrow[u,"\Lambda_R"'] & \mu'_{2}(R) \arrow[r] \arrow[u] & 1.
\end{tikzcd}
    \end{center}
 Let $\alpha\in\K_{2}(\D_{r},R)$ such that $\Gamma_{R}(\alpha)=1$. Since the square commutes, $\Lambda_R(j(\alpha))=i(\Gamma_{R}(\alpha))=1$. Note that the middle vertical row $\Lambda_R$ is injective by Theorem $\ref{thm_stb_dr_panin}$. Thus $\alpha=1$, and hence $\Gamma_{R}$ is injective.
\end{proof}
\section{Stability Theorems over Laurent Polynomial rings}\label{sec_lpr}
\begin{lem}(cf. \cite[Theorem 2.1]{Ri80})
    Let $R$ be a Noetherian ring, and $S$ denote the multiplicative closed subset of special monic polynomials in $R[X]$. If $Q$ is a maximal ideal of $R[X,X^{-1}]$ with $\hh(Q)>\dim(R)$, then $Q\cap S\neq \phi$. In particular, $\dim(S^{-1}R[X,X^{-1}])\leq \dim(R)$.\label{lem_dim_lpr_spmonic}
\end{lem}
\begin{proof}
 Since $Q$ is a maximal ideal of $R[X,X^{-1}]$ with $\hh(Q)>\dim(R)$, it is easy to observe that there is a maximal ideal $P$ of $R[X]$ with $\hh(P)>\dim(R)$ such that $P_{X}=Q$. Since $P$ is maximal and $X\notin P$, $P+(X)=R[X]$. Hence, there exist $f\in R[X]$ and $p\in P$ such that $p+fX=1$. Now, evaluating $X=0$, we get $p(0)=1$. By \cite[Lemma 6.2]{Su77b}, there exists a monic polynomial $g\in P$. Set $g(0)=a$. If $\deg(g)>\deg(p)$, then we choose $h=g-(a-1)p$. If $\deg(g)\leq \deg(p)$, we can consider large $n$ such that $\deg(g^n)>\deg(p)$; then we choose $h=g^n-(a^n-1)p$. In both cases, we have a special monic polynomial $h$, and hence $h\in Q\cap S$.
\end{proof}
\begin{thm}\label{thm_surj_stab_Laur_poly}
     Let $R$ be a Noetherian domain. Let $\Phi'\subset \Phi$ be embeddings of root systems as in the statement of Theorem $\ref{thm_surj_stab_poly}$. Assume that $\K_{1}(\Phi,R[X])$ is trivial. Then the canonical homomorphism $\K_{2}(\Phi',R[X,X^{-1}])\to \K_{2}(\Phi,R[X,X^{-1}])$ is surjective.
\end{thm}
\begin{proof}
    Let $S$ denote the multiplicative closed subset of special monic polynomials in $R[X]$. From Theorem $\ref{thm_sr_bass}$ and Lemma $\ref{lem_dim_lpr_spmonic}$, we have $$\sr(S^{-1}R[X,X^{-1}])\leq\dim(S^{-1}R[X,X^{-1}])+1\leq\dim(R)+1.$$ For (vii), both the ring $R$ and $S^{-1}R[X,X^{-1}]$ have dimension $\leq 1$. Now, applying \cite[Theorem 5.9]{Si25} (for base case) to $R$ and $S^{-1}R[X,X^{-1}]$, we conclude that the maps $\K_{2}(\A_4,R)\to \K_{2}(\Psi,R)$ and $\K_{2}(\A_4,R[X,X^{-1}]_f)\to \K_{2}(\Psi,R[X,X^{-1}]_f)$ are surjective for some $f\in S$, where $\Psi$ is a root system of type $\A_{\geq 4},~\D_{\geq 5}$ or $\E_{6,7,8}$. \par For (i)--(vi), using a similar argument as in Step 1 and Step 2 of Theorem $\ref{thm_surj_stab_poly}$, we conclude that $\K_{2}(\Phi',R)\to \K_{2}(\Phi,R)$ and $\K_{2}(\Phi',R[X,X^{-1}]_f)\to \K_{2}(\Phi,R[X,X^{-1}]_f)$ are surjective for some $f\in S$. Then, for all cases, the result follows from Theorem $\ref{thm_inj_surj_Laur_poly_spmonic}(2)$.
\end{proof}
We now mention below some classes of rings $R$ for which $\K_{1}(\Phi,R[X])$ is trivial. It is well known that $\K_{1}(\Phi,R)=1$ for any local ring or Euclidean domain $R$ (cf. \cite{Ab69,Ab83}). In \cite[Theorem 3.5]{Ab83}, Abe proved that $\K_{1}(\Phi,k[X_1,\ldots,X_n])=1$, where $k$ is a field. Later, Stavrova generalized this result to reductive groups over various rings (see \cite{Sta20,Sta22,Sta24} for more details). Before giving examples, we first recall the following definitions from \cite[Section 1]{Sta20} and \cite[Chapter 4.6]{Lambook06}, respectively.
\begin{defn}
    A ring $R$ is said to be a Dedekind domain of arithmetic type if $R=\mathcal{O}_{k,S}$, is the ring of $S$-integers of a global field $k$ with respect to a finite non-empty set $S$ of primes containing all archimedean primes. Typical examples include $\Z$, function field $\mathbb{F}_{p^{m}}[t]$, and localization of rings of integers at finitely many primes, such as $\Z_{(p)}$, where $p$ is a prime.
\end{defn}
\begin{defn}
     A PID $R$ is said to be a special PID if $\SL_{n}(R)=\E_{n}(R)$ for all $n\geq 2$. The class of special PIDs contains Euclidean domains and localizations of 2-dimensional regular local rings at a regular parameter, including the rings $A(x)$ for a discrete valuation ring $A$. Any localization of a special PID is also a special PID.
\end{defn}
\begin{exms}\hfill\break
    \begin{enumerate}
        \item If $R$ is a Dedekind domain of arithmetic type, then $\K_{1}(\Phi,R[X])=1$ (cf. \cite[Corollary 1.2]{Sta20}). Recall that $\rank(\Phi)\neq 1$, otherwise $\K_{1}(\A_1,\Z[X])$ is not necessarily equal to $1$ since $\SL_2(\Z[X])\neq\E_2(\Z[X])$. Note that, if $R$ is a Dedekind domain (not necessarily of arithmetic type), then $\K_{1}(\A_{2},R)\neq 1$ (see \cite[page 57]{Lambook06}).
        \item If $R$ is a Dedikind ring, then $\K_{1}(\Phi,R[X])\cong \K_{1}(\Phi,R)$ (cf. \cite[Theorem 2.6]{Sta24}). In addition if $R$ is a local ring, then $\K_{1}(\Phi,R[X])=1$. In particular, if $R=\mathbb{F}_{2}~\text{or}~\mathbb{F}_{2}[t]$ (i.e., no non-trivial units), then $\K_{1}(\Phi,R[X])=1$.
        \item If $R$ is a special PID, then $\K_{1}(\Phi,R[X])=1$ (cf. \cite[Theorem 1.2]{Sta24}).  
        \item If  $R$ is a regular ring containing a field, then $\K_{1}(\Phi,R[X])\cong \K_{1}(\Phi,R)$ (cf. \cite[Theorem 1.1]{Sta22}). In addition if $R$ is a local ring, then $\K_{1}(\Phi,R[X])=1$.
        \item If $R$ is an equicharacteristic regular complete local ring or a regular local ring of essentially finite type over a field, then $\K_{1}(\Phi,R[X])=1$ (cf. \cite[Corollary 3.9]{Ab83})).
        \item If $R$ is a Noetherian ring geometrically regular over a Dedekind ring with perfect residue fields, then $\K_{1}(\Phi,R[X])\cong \K_{1}(\Phi,R)$ (cf. \cite[Theorem 2.6]{Sta24}). In addition if $R$ is a local ring, then $\K_{1}(\Phi,R[X])=1$.
    \end{enumerate}
\end{exms}
\begin{cor}
    Let $R$ be a Noetherian domain, and let $\Phi'\subset \Phi$ be embeddings of root systems as in the statement of Theorem $\ref{thm_surj_stab_poly}$. Assume that $R$ is one of the following:
    \begin{enumerate}
        \item[$(i)$] a Dedekind domain of arithmetic type,
        \item[$(ii)$] a special PID,
        \item[$(iii)$] a regular local ring containing a field,
        \item[$(iv)$] an equicharacteristic regular complete local ring,
        \item[$(v)$] a regular local ring of essentially finite type over a field,
        \item[$(vi)$] a Noetherian ring geometrically regular over a Dedekind ring with perfect residue fields,
        \item[$(vii)$] a Dedekind local ring. 
    \end{enumerate}
    Then the canonical homomorphism $\K_{2}(\Phi',R[X,X^{-1}])\to \K_{2}(\Phi,R[X,X^{-1}])$ is surjective.
\end{cor}
At present, we do not know whether Theorem $\ref{thm_surj_stab_Laur_poly}$ holds without assuming the triviality of $\K_{1}(\Phi,R[X])$. In particular, we would like to know the following. 
\begin{ques}\label{ques_LPR}
   Let $R$ be a Noetherian ring and $A=R[X,\ldots,X_n,Y_{1}^{\pm1},\ldots,Y_{m}^{\pm1}]$. Let $\Phi'\subset \Phi$ be embeddings of root systems as in the statement of Theorem $\ref{thm_surj_stab_poly}$. Does the canonical homomorphism $\K_{2}(\Phi',A)\to \K_{2}(\Phi,A)$ surjective?  
\end{ques}
\begin{thm}\label{thm_inj_stab_Laur_poly}
    Let $R$ be a Noetherian domain and $A=R[X,\ldots,X_n,Y_{1}^{\pm1},\ldots,Y_{m}^{\pm1}]$. Then the canonical homomorphism $\K_{2}(\A_{r},A)\to \K_{2}(\A_{r+1},A)$ is injective for $r\geq\max$ $\{4,\dim(R)+2\}$.
\end{thm}
\begin{proof}
For $m=0$, the result follows from Theorem $\ref{thm_inj_stab_poly}$. So we assume $m\geq 1$. Let $S$ denote the multiplicative closed subset of special monic polynomials in $R[Y_m]$. Then $$S^{-1}A=(S^{-1}R[Y_m,Y_m^{-1}])[X_1,\ldots,X_n,Y_{1}^{\pm1},\ldots,Y_{m-1}^{\pm1}],$$ where $\dim(S^{-1}R[Y_m,Y_m^{-1}])\leq \dim(R)$. By induction step, the map $\K_{2}(\A_{r},S^{-1}A)\to \K_{2}(\A_{r+1},S^{-1}A)$ is injective, and hence there exists $f\in S$ such that $\K_{2}(\A_{r},A_f)\to \K_{2}(\A_{r+1},A_f)$ is injective. Since $A_f=R[X_1,\ldots,X_n,Y_{1}^{\pm1},\ldots,Y_{m-1}^{\pm1}][Y_m,Y_m^{-1}]_f$, the result follows from Theorem $\ref{thm_inj_surj_Laur_poly_spmonic}(1)$.
\end{proof}
The following result is an unstable $\K_2$-homotopy invariance property for Laurent polynomial rings.
\begin{thm}\label{thm_hom_inv_Laur_poly}
     Let $R$ be a Dedekind domain or field, having no non-trivial units. Let $\Psi$ be a root system of type $\A_{\geq 4},~\D_{\geq 5}$ or $\E_{6,7,8}$.  Assume that $\K_{1}(\Psi,R)$ is trivial. Then $\K_{2}(\Psi,R[X,X^{-1}])\cong \K_{2}(\Psi,R)$.
\end{thm}
\begin{proof}
 Recall that $\K_{2}(R)=\underset{\to r}{\lim}~\K_{2}(\A_r,R)$ and $\K_{1}(R)\cong \SK_{1}(R) (=\underset{\to r}{\lim}~\K_{1}(\A_r,R)=1)$, since the unit of $R$ is trivial. Consider the following commutative diagram:
    \begin{center}
        \begin{tikzcd}
\K_{2}(\A_4,R) \arrow[r,"\cong"] \arrow[d, "\xi"] & \K_{2}(R) \arrow[d,"\cong"] \\
\K_{2}(\A_4,R[X,X^{-1}]) \arrow[r,"\cong"]                  & \K_{2}(R[X,X^{-1}]).          
\end{tikzcd}
\end{center}
From \cite[Theorem 2.6]{Sta24}, we have $\K_{1}(\Psi,R[X])=1$. The two horizontal arrows in the above diagram are isomorphisms by combining Theorems $\ref{thm_surj_stab_poly},\ref{thm_inj_stab_poly},\ref{thm_surj_stab_Laur_poly}$ and $\ref{thm_inj_stab_Laur_poly}$, as we have a sequence of maps $\K_{2}(\A_r,-)\to\K_{2}(\A_{r+1},-)$ for $r\geq 4$. The right vertical arrow is an isomorphism by stable $\K_2$-invariance (cf. \cite[Section V, Theorem 6.3]{Weibel-3}). Therefore, the map $\xi$ is an isomorphism. Now consider another commutative diagram:
    \begin{center}
        \begin{tikzcd}
\K_{2}(\A_4,R) \arrow[r] \arrow[rrd, "\xi"'] & \K_{2}(\Psi,R) \arrow[r, "\zeta"] & \K_{2}(\Psi,R[X,X^{-1}])           \\
                               &                      & \K_{2}(\A_4,R[X,X^{-1}]). \arrow[u,two heads]
\end{tikzcd}
    \end{center}
    By Theorem $\ref{thm_surj_stab_Laur_poly}$, the vertical map in the above diagram is surjective, and hence the map $\zeta$ is surjective. It is easy to verify that the map $\zeta$ is split by the obvious homomorphism given by evaluation at $X=1$, and therefore it is injective. This completes the proof.
\end{proof}
In the rest of this section, we provide a similar injectivity result using localization by inverting monic polynomials in $R[X]$ and $R[X^{-1}]$, instead of a special monic polynomial.
\begin{prop}\label{prop_inj_Laur_poly_monic2}
    Suppose $r\geq 4$. Let $R$ be a Noetherian domain and $A=R[X,X^{-1}]$. Suppose $f_{+}(X)\in A[X]$ (resp. $f_{-}(X^{-1})\in A[X^{-1}]$) is a monic polynomial. If both homomorphisms $\K_{2}(\A_{r},A_{f_{\pm}})\to \K_{2}(\A_{r+1},A_{f_{\pm}})$ are injective, then the homomorphism $\K_{2}(\A_{r},A)\to \K_{2}(\A_{r+1},A)$ is injective.
\end{prop}
\begin{proof}
Write $f_{-}(X^{-1})=X^{-h}+c_1X^{-h+1}+\cdots+c_h$ and $g=X^{h}f_{-}(X^{-1})=c_hX^{h}+\cdots+c_1X+1$. Then $X$ and $g$ are comaximal elements of $R[X]$. Consider the following Cartesian square (I): 
   \begin{center}
       \begin{tikzcd}
R[X] \arrow[r,"\lambda_{X}"] \arrow[d,"\lambda_{g}"] & R[X,X^{-1}] \arrow[r,"\lambda_{f_{+}}"] \arrow[d,"\lambda_{f_{-}}"] & R[X,X^{-1}]_{f_{+}} \\
R[X]_{g} \arrow[r,"\lambda_{X}"]           & R[X,X^{-1}]_{g}=R[X,X^{-1}]_{f_{-}}   
\arrow[phantom, from=1-1, to=2-2, "\text{I}" description]
\end{tikzcd}
\end{center}
    Let $\alpha\in \K_{2}(\A_{r},A)$ such that $l_{A}(\alpha)=1$ in $\K_{2}(\A_{r+1},A)$. Then $\lambda_{f_{\pm}}^{\A_{r+1}}(l_{A}(\alpha))=1$, and hence $l_{A_{f_{\pm}}}(\lambda_{f_{\pm}}^{\A_{r}}(\alpha))=1$ due the commutative diagram. Since $l_{A_{f_{-}}}$ is injective, $\lambda_{f_{-}}^{\A_{r}}(\alpha)=1=\lambda_{X}^{\A_{r}}(1)$. Applying Proposition $\ref{prop_excision_k2}(2)$, we get $\beta\in\K_{2}(\A_{r},R[X])$ such that $\lambda_{X}^{\A_{r}}(\beta)=\alpha$ and $\lambda_{g}^{\A_{r}}(\beta)=1$. Then $l_{A_{f_{+}}}(\lambda_{Xf_{+}}^{\A_{r}}(\beta))=l_{A_{f_{+}}}(\lambda_{f_{+}}^{\A_{r}}(\alpha))=1$ and since $l_{A_{f_{+}}}$ is injective, $\lambda_{Xf_{+}}^{\A_{r}}(\beta)=1$. Since $Xf_{+}$ is a monic polynomial in $R[X]$, by \cite[Corollary 5.2]{Tu83}, $\beta=1$, and hence $\alpha=1$. 
\end{proof}
We now prove Theorem $\ref{thm_inj_stab_Laur_poly}$ in a similar manner using Proposition $\ref{prop_inj_Laur_poly_monic2}$ instead of Theorem $\ref{thm_inj_surj_Laur_poly_spmonic}$.\par
\textbf{An alternative proof of Theorem $\ref{thm_inj_stab_Laur_poly}$.}
     We denote by $S_{+}$ and $S_{-}$ the multiplicative closed subsets of monic polynomials in $R[Y]$ and $R[Y^{-1}]$ respectively. Then $\dim(S_{+}^{-1}R[Y,Y^{-1}])=\dim(S_{+}^{-1}R[Y])\leq\dim(R)$. Hence $\sr(S_{+}^{-1}R[Y,Y^{-1}])\leq\dim(R)+1$. Since $r\geq\max\{4,\sr(S_{+}^{-1}R[Y,Y^{-1}])+1\}$, by Theorem $\ref{thm_sutu_inj}$ the homomorphism $$\K_{2}(\A_{r},S_{+}^{-1}R[Y,Y^{-1}])\to \K_{2}(\A_{r+1},S_{+}^{-1}R[Y,Y^{-1}])$$ is injective. Similarly, the homomorphism $$\K_{2}(\A_{r},S_{-}^{-1}R[Y,Y^{-1}])\to \K_{2}(\A_{r+1},S_{-}^{-1}R[Y,Y^{-1}])$$ is also injective. Then use Proposition $\ref{prop_inj_Laur_poly_monic2}$ and apply induction on $m$.\hfill \qedsymbol \par
 %Similarly, we have $\sr(S_{-}^{-1}R[X,X^{-1}])\leq\dim(R)+1$.
We conclude this section with the following question, which naturally suggests a possible generalization of Question $\ref{ques_LPR}$ to the setting of monoid algebras.
\begin{ques}
     Let $R$ be a Noetherian ring and $M$ be a commutative cancellative monoid. Let $\Phi'\subset \Phi$ be embeddings of an irreducible root system (not necessarily simply-laced). Does the canonical homomorphism $\K_{2}(\Phi',R[M])\to \K_{2}(\Phi,R[M])$ surjective? 
\end{ques}
\section{Concluding Remarks} In this section, we will discuss two important consequences of the stability results suggested by the reviewer. \par
  $\blacktriangleright$ One can ask whether surjective $\K_2$-stability results can provide an alternative route to the centrality of $\K_2(\Phi,R)$ in $\St(\Phi,R)$? \par
   First, we will discuss this in the linear case, i.e., $\Phi$ is of type $\A_r$. If the map $\K_2(\A_{r-1},R)\to \K_2(\A_{r},R)$ is surjective for $r\geq 3$, then by \cite[1.4.5]{HaMe89}, $\K_2(\A_{r},R)$ is central in $\St(\A_{r},R)$. However, the converse is not true. van der Kallen \cite{Ka77} proved that $\K_2(\A_{r},R)$ is central in $\St(\A_{r},R)$ for $r\geq 3$, but the map $\K_2(\A_{r-1},R)\to \K_2(\A_{r},R)$ may not be surjective ( for example, see \cite[Section 7]{Ka76}). \par Now, for a general root system, the following result is due to Stein \cite[Theorem 5.1]{St71}. 
   \begin{thm}\label{thm_stein_central}
     Let $\Delta$ denote a set of fundamental roots in a root system $\Phi$. For $\alpha\in \Delta$, if $\Phi'\subset\Phi$ is the subsystem generated by $\Delta\setminus \alpha$, then \[\K_2(\Phi,R)\cap \Img(\St(\Phi',R)\to \St(\Phi,R))\] is central in $\St(\Phi,R)$.
   \end{thm}
Now, if the map $\K_2(\Phi',R)\to \K_2(\Phi,R)$ is surjective for a subsystem $\Phi'$ obtained by deleting a suitable simple root from $\Phi$ (i.e., deleting one node from the Dynkin diagram of $\Phi$), then it is easy to see that $\K_2(\Phi,R)\subseteq \Img(\St(\Phi',R)\to \St(\Phi,R))$, and hence, by Theorem $\ref{thm_stein_central}$, \[\K_2(\Phi,R)\cap \Img(\St(\Phi',R)\to \St(\Phi,R))=\K_2(\Phi,R)\] is central in $\St(\Phi,R)$. However, the centrality of $\K_2(\Phi,R)$ for any root system $\Phi$ of rank $\geq 3$ is already known due to the work of various authors (see \cite{LaSiVo24}). Therefore, this answers the given question. It is observed that all ``dasharrow" embeddings of root systems in Figure $\ref{root_embedding}$ satisfy the stated hypothesis.\par

$\blacktriangleright$ Considering stability results, can one speculate on the important question of when the unstable $\K_2$-group for a root system $\Phi$ is generated by Steinberg symbols, e.g., over Dedekind domains such as rings of integers in number fields?\par
\begin{defn}
  For $\alpha\in \Phi$ and $s\in R$, define certain elements $w_{\alpha}(s)$, $h_{\alpha}(s)$ of $\St(\Phi,R)$ as \[w_{\alpha}(s)=x_{\alpha}(s)\cdot x_{-\alpha}(-s^{-1})\cdot x_{\alpha}(s),~~~  h_{\alpha}(s)=w_{\alpha}(s)\cdot w_{\alpha}(-1).\]
  The Steinberg symbols are defined for arbitrary $s,t\in R^{\times}$ as follows: $$\{s,t\}_{\alpha}=h_{\alpha}(st)\cdot h^{-1}_{\alpha}(s)\cdot h_{\alpha}(t).$$ In turn, the Dennis–Stein symbols are defined for arbitrary $a,b\in R$ satisfies $1+ab\in R^{\times}$ as follows: $$\left<a,b\right>_{\alpha}=x_{-\alpha}\left(\frac{-b}{1+ab}\right)\cdot x_{\alpha}(a)\cdot x_{-\alpha}(b)\cdot x_{\alpha}\left(\frac{-a}{1+ab}\right)\cdot h^{-1}_{\alpha}(1+ab).$$
\end{defn}  
  Since every Steinberg symbol can be expressed as a Dennis–Stein symbol via the identity $\{s,t\}_{\alpha}=\left<-s,\frac{1-t}{s}\right>_{\alpha}$ for $s,t\in R^{\times}$, it follows that if $\K_2(\Phi,R)$ is generated by Steinberg symbols then it is also generated by Dennis–Stein symbols.\par
  Several examples and counterexamples addressing this question have been provided by various authors \cite{AbMo88,BrHu82,ChYu03,Hu83,Hu99,Ka77,Ka77a,Ma69,Mu92,St73}. We now present a few results to further illustrate this question.
  \begin{itemize}
        \item It is known that Milnor's $\K_2$-group $\K_2(R)$ is generated by Steinberg symbols when the underlying ring $R$ is a field, a division ring, a local ring, or even a semilocal ring \cite[Theorem 5.10.5]{Weibel-3}. Using Theorem $\ref{thm_sr_stein}$ and Theorem $\ref{thm_sutu_inj}$, the unstable $\K_2$-group $\K_2(\A_r,R)$ is generated by Steinberg symbols for $r\geq \sr(R)+1$.
      \item Let $\mathcal{O}_K$ denote the ring of integers of a number field $K/\mathbb{Q}$. Browkin-Hurrelbrink \cite{BrHu82} showed that $\K_2(\mathcal{O}_{\mathbb{Q}(\sqrt{d})})$ is generated by Steinberg symbols if and only if $d=2,5,13$. Now, using van der Kallen stability result \cite[Theorem 1]{Ka81}, one concludes that $\K_2(\A_{\geq 2},\mathcal{O}_{\mathbb{Q}(\sqrt{d})})$ is generated by Steinberg symbols if and only if $d=2,5,13$.
      \item Mulders \cite[Proposition 6.9]{Mu92} showed that $\K_2(\mathcal{O}_{\mathbb{Q}(\sqrt{d})})$ is generated by Dennis–Stein symbols provided that $d\neq \pm 1$ (mod 8) and $d\neq 6$ (mod 9). Now, using \cite[Theorem 1]{Ka81}, one shows that $\K_2(\A_{\geq 2},\mathcal{O}_{\mathbb{Q}(\sqrt{d})})$ is generated by Dennis–Stein symbols provided that $d\neq \pm 1$ (mod 8) and $d\neq 6$ (mod 9).
       \item Hutchinson constructs an explicit biquadratic field extension $K=\mathbb{Q}(\sqrt{-2d},\sqrt{f})$ (under certain conditions on $d$ and $f$, e.g., $K=\mathbb{Q}(\sqrt{-34},\sqrt{-206}))$ and proves that $\K_2(\mathcal{O}_K)$ is not generated by Dennis–Stein symbols (for more details, see \cite[Lemma 4.2]{Hu99}). Consequently, in this case $\K_2(\A_{\geq 2},\mathcal{O}_K)$ is also not generated by Dennis–Stein symbols and hence not generated by Steinberg symbols.
       \item $\K_2(\A_1,\mathcal{O}_{\mathbb{Q}(\sqrt{d})}[\frac{1}{p}])$ ($d$ depends on the prime number $p$) cannot be generated by Steinberg symbols; see \cite[Theorem 1]{ChYu03}.
       \item Hurrelbrink \cite[Theorem 1]{Hu83} showed that $\K_2(A_{\geq 2},\mathbb{Z}[\frac{1}{p_1},\ldots,\frac{1}{p_n}])$ is generated by Steinberg symbols and Dennis–Stein symbols where $p_i$ are prime numbers.
     % \item In \cite{BeGa04}, the authors provide tables illustrating when $K_2$ is generated by Steinberg symbols over the ring of integers of a quadratic field extension.
   \end{itemize}
Therefore, in view of the above discussion, one may ask when the group $\K_2(\Phi,\mathcal{O}_K)$ generated by Steinberg symbols for a root system $\Phi$ not of type $\A$?\par
 For the general root system $\Phi$, the group $\K_2(\Phi,R)$ is known to be generated by Steinberg symbols in the following cases:
 \begin{itemize}
     \item Let $k$ be a field. The case $R=k$ is due to Matsumoto \cite[Corollaire 5.11]{Ma69}, $R=k[X]$ to Rehmann \cite[Corollary 6]{Re75} and $R=k[X,X^{-1}]$ to Abe-Morita \cite[Corollary 6]{AbMo88}.
     \item The case where $R$ is semilocal with at most one residue field isomorphic to $\mathbb{F}_2$ and $\rank(\Phi)\geq 2$ is due to Stein \cite[Corollary 2.14]{St73}.
     \item $R=\mathbb{Z}[\frac{1}{p_1},\ldots,\frac{1}{p_n}]$ under certain conditions on $p_i$, see \cite[Corollary 7]{AbMo88}.
     \item The case where $R$ is U-irreducible is due to van der Kallen \cite[Theorem 3.4, 3.7]{Ka77a}.
 \end{itemize}
\textbf{Acknowledgments} I sincerely thank the anonymous referee for their careful reading and helpful remarks. I am deeply grateful to Husney Parvez Sarwar for introducing me to the field of $\K$-theory and to Manoj Kumar Keshari for listening to my lectures on this work and offering numerous corrections. I would also like to thank Sergey Sinchuk for his detailed comments on the first draft of the paper, as well as Shripad M. Garge and Anastasia Stavrova for clarifying my doubts. This work was supported by the IIT Bombay Postdoctoral Research Fellowship.

%\textbf{\large Declarations} \\

%\textbf{Conflict of Interest} The authors have no conflict of interest to declare that are relevant to this article.\\ \\ 
%\textbf{\large Data availability} \\

% Data sharing not applicable to this article as no datasets were generated or analysed during the current study.
 %No data was used for the research described in the article.

%Ann. K-Theory (Marco Schilithing),Journal of Algebra, JPAA, C. R. Acad. Sci. Paris, Algebra and representation theory, Bulletin des Sciences Mathématiques Indranil Biswas, documenta mathematices ulf Rehaman
\end{document}